%% file: main.tex
\documentclass[11pt,a4paper]{article}

\title{Efficient prediction of static and dynamical responses of functional graded beams using sparse multiscale patches}

\author{%
Thien Tran-Duc\thanks{%
Mathematical Sciences, University of Adelaide,  South Australia.
\protect\url{mailto:thien.tran-duc@adelaide.edu.au},
\protect\url{https://orcid.org/0000-0002-2004-5156}}
\and 
J. E. Bunder\thanks{%
UniSA STEM, University of South Australia, Australia.
\protect\url{mailto:Judy.Bunder@unisa.edu.au},
\protect\url{https://orcid.org/0000-0001-5355-2288}}
\and 
A. J. Roberts\thanks{%
University of Adelaide,  South Australia.
protect\url{mailto:ProfAJRoberts@protonmail.com},
\protect\url{https://orcid.org/0000-0001-8930-1552}}
}

\input{preamble}

\IfFileExists{ajrxxx.sty}%
   {\let\inPlot\input}%
   {\let\inPlot\includegraphics}

\begin{document}

\maketitle

\begin{abstract}
We develop a multiscale patch scheme for studying the system level characteristics of heterogeneous functional graded beams.
The algorithm computes the detailed beam dynamics on the microscale, but only in small \emph{patches} of the beam domain, and then applies symmetry-preserving interpolation to these patches to accurately predict the macroscale behaviour.
To validate the algorithm, two examples of functionally graded beams are investigated, namely cross-sectionally graded and axially graded.
Gradient patterns are defined via volume fractions of aluminium and silicon carbine either over the beam's cross section or along its axial direction. 
In these examples the multiscale patch scheme only computes over a fraction of the beam's full-domain. 
Beam deflection and natural frequencies from the patch computations agree very well with existing experimental data and the full-domain computations.
The algorithm is stable and robust, with errors consistently small and reliably reducible by increasing the number of patches.
The reduction in the spatial domain of computation substantially improves the computational efficiency, with the computational time reducing by a factor of up to~\(17\) when the patches cover~$27\%$ of the beam.
\end{abstract}

\tableofcontents

\section{Introduction}


Functional graded materials (\textsc{fgm}s) are constructed from multiple material components in a certain gradient pattern in one or more material dimensions.
Such materials are typically designed to have particular characteristics, such as extreme mechanical strengths or physical and chemical properties for specific applications \citep[e.g.,][]{Torquato2010, Somnic2022}. 
As physical experimentation is time-consuming and expensive, it is becoming increasingly common for numerical simulations to be the first step in designing new \textsc{fgm}s.
However, predicting the macroscale behaviour is difficult as simulations need to resolve the detailed microscale structure.  
The Equation-Free Patch Scheme \citep[e.g.,][]{Kevrekidis2003, Kevrekidis2004, Samaey2005, Samaey2006} is a computationally efficient,  controllably accurate, method for macroscale prediction of such multiscale systems.

\begin{figure}\centering
\caption{\label{FAxiallyGradedBeamSim1}six snapshots in time of the bending and twisting of a 3D heterogeneous, functionally graded \Al/\Si\C\ beam (non-dimensional) following the initial bent state shown at time \(t=0\).  Physically, the times here are roughly in milliseconds, and the colour represents the elastic stress component~\(\sigma^{xx}\).  
Computation is performed only within the nine shown small patches \text{in 3D space.}}
\includegraphics{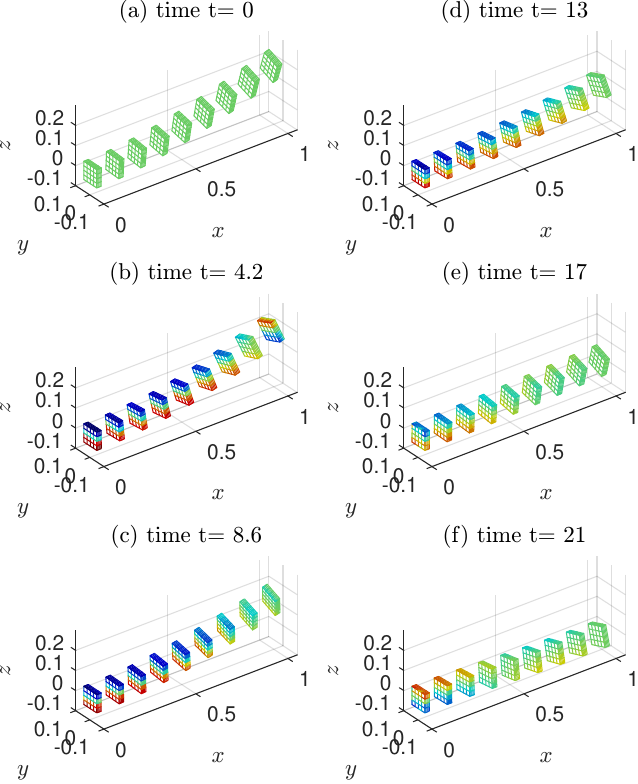}
\end{figure}%
For example, \cref{FAxiallyGradedBeamSim1} shows six snapshots in time of the bending and twisting in time of a 3D heterogeneous functionally graded beam (\cref{Sec:AxialGB}) with one end fixed at \(x=0\) and one end free near \(x=1\)\,.
The simulation starts at time \(t=0\) (top-left) with some bending but an assumed zero stress \emph{within} each patch.  
Then for times \(t>0\) coupling between the patches (\cref{sec:coupling}) provides the appropriate stress on the left and right faces of the patches (the faces nearly perpendicular to the \(x\)-direction, roughly in the \(yz\)-plane).
These face-stresses (represented by colour in \cref{FAxiallyGradedBeamSim1}) rapidly propagate across the microscale of each patch to generate the bending moments apparent in subsequent times of \cref{FAxiallyGradedBeamSim1}, bending moments that accelerate the beam downwards.
\cref{FAxiallyGradedBeamSim1} shows about a quarter of one period of the bending vibration.
The initial condition also has some twist, and the six snapshots show about~\(2.5\) periods of the simultaneous twisting of the beam.
Such accurate time dynamic simulation and prediction is difficult for most multiscale schemes, but is straightforward with the developed and proven equation-free patch scheme that we further develop herein \citep[e.g,][]{Maclean2020a, Roberts2019b}.
%

Functionally graded materials (\textsc{fgm}s) are often classified into two types: sandwich layered patterns \citep{Gao2020, Nian2023}; and continuous gradient patterns \citep{Kumar2022, Pasha2022, Ganta2022}.
Of the two, sandwich layered materials have received more attention because they are easier to produce as they are manufactured by combining layers of different types of materials \citep[e.g.][]{Gao2020,Nian2023,Peng2021}.
Such sandwich structures are designed to have useful properties, such as, light weight, high stiffness-to-weight ratio, excellent energy absorption efficiency, and advanced thermal or mechanical properties \citep{Chen2024, Li2019, Garg2022, Wei2020}.  
\cref{Sec:CrSecGB} applies the patch scheme to a sandwich layered beam.
The beams consist of multiple layers of mixtures of aluminium and silicon carbide, arranged with increasing volume fractions of silicon carbide.
The patch scheme is validated via comparison with the experimental data of \cite{Kapu2008}.

The other class of \textsc{fgm}s, continuous gradient patterns,  are made from mixtures of two or more materials such that physical properties  vary continuously across one or more dimensions.
These materials are often made via 3D printing processes so that the material blending is well controlled.  
These  \textsc{fgm}s are commonly metal-ceramic mixtures that inherit the high rigidity, thermal stability and chemical resistance of the ceramic component, and the highly tensile property of the metal component.
Such properties are ideal for high tensile stress or high temperature environments, such as in space or bio-medical  applications \citep{Lin2009,Pompe2003,Sola2016}.
\cref{Sec:AxialGB} applies the patch scheme to a \textsc{fgm} with a continuous gradient pattern, and leads to the simulation \text{illustrated in \cref{FAxiallyGradedBeamSim1}.}

Homogenisation is a widely-used technique in modelling heterogeneous elasticity \cite[e.g.,][]{Fish2021}.
Homogenisation techniques variously invoke  
algebraic analysis of representative volume elements (\rve) \cite[e.g.,][]{RamirezTorres2018, Pau2022},
numerical precomputation of a library of canonical deformations \cite[e.g.,][]{Raju2021, Schneider2021}, or
computational on-the-fly closures \cite[e.g.,][]{Abdulle2012, Kevrekidis2004}.
The equation-free patch scheme developed herein (\cref{Sec:Patchscheme}), and seen in \cref{FAxiallyGradedBeamSim1}, provides an efficient and accurate computational homogenisation \citep[e.g.,][]{Kevrekidis2003, Kevrekidis2004, Samaey2005, Samaey2006}.
A novel distinguishing feature of the patch scheme is that the underlying theory \cite[e.g.,][]{Roberts06d, Roberts2011a, Bunder2020a} applies rigorously at the \emph{finite} scale separation of real physical scenarios.
This support at finite scale separation is important, as in many practical beam modelling scenarios the scale separation parameter may be as large as \(0.1\)--\(0.3\).
Further, theory for heterogeneous systems \cite[e.g.,][]{Bunder2020a} shows that the patch scheme has controllable accuracy at finite scale separation, even up to an almost unprecedented eight significant digits in the macroscale predictions. 
The multiscale patch scheme achieves high computational efficiency by only computing the microscale complex physical phenomena in small, sparse, patches of space \citep[or more generally in small, sparse, patches of \emph{space-time},][]{Kevrekidis2004}.
Patches are analogous to the \rve{}s of other methodologies.
The patches are positioned in space to resolve the desired macroscale phenomena. 
\cref{sec:coupling} discusses how patches are coupled via inter-patch interpolation which determines field values on the 3D patch faces.
The patch coupling conditions are distinct and unrelated to the physical boundary conditions on the macroscale domain of the system. 
In all cases, the patch faces \emph{not} on a physical boundary have values set by the inter-patch interpolation of \cref{sec:coupling}.
Then, the microscale simulations within each patch automatically and \emph{implicitly} provide a `high-order' accurate homogenised closure on the scale of the macroscale grid of patches.
Thus the patch scheme provides a computationally efficient homogenisation.

In a recent review of multiscale modelling in materials, \cite{Fish2021} [p.774] comment that pursuing ``a multiscale approach involves a trade-off between increased model fidelity with the added complexity, and corresponding reduction in precision and increase in uncertainty''.
A significant property of the multiscale patch scheme is that there is no need for such a trade-off.
Further, \cite{Fish2021} [p.781] assert that ``Progress in this and similar methods requires research on accurate and efficient lifting operators for specific applications.''
Theory underlying the patch scheme proves \citep{Bunder2020a}, and the elastic beam application developed herein confirms, that the patch coupling of \cref{sec:coupling} implicitly provides such ``accurate and efficient lifting'' \text{in general circumstances.}

The equation-free patch scheme is designed to wrap around any given microscale system to then non-intrusively empower efficient and accurate macroscale prediction.
A freely available Equation-Free Toolbox for \matlab\ or Octave \citep{Maclean2020a, Roberts2019b} invokes a user-defined function of the microscale physics, and applies the patch scheme as a `black-box'.
A user also specifies the macroscale geometry and boundary conditions of interest in space-time.
The patch scheme is designed for, and has been tested on, nonlinear systems; however, herein we restrict attention to linear viscoelasticity.  
The macroscale results of \cref{Sec:CrSecGB,Sec:AxialGB} are all obtained by  applying functions of the Equation-Free Toolbox on the given microscale \text{viscoelasticity of \cref{sec:discretisation}.}

To summarise this introduction, the power of the patch scheme for computational homogenisation in elasticity is that the scheme is computationally efficient, with controllable accuracy---potentially controllable down to round-off errors, stable for dynamic simulations of wave problems,  non-intrusive in that it just `wraps around' a user's code, can handle nonlinear systems, does not require a variational principle,  has proven accuracy at the finite scale-separation of real physics,  is proven to implicitly provide a correct macroscale closure,  requires no expensive preprocessing,  there are no assumed boundary conditions on \rve{}s---instead there are proven coupling conditions, there is no oversampling regions, no buffer regions, no action regions, no need to guess specific fast\slash slow variables, and  no need for a trade-off between fidelity and complexity.
The many other extant multiscale approaches each have some of these desirable properties. 
The novelty of our patch scheme is that no other approach shares all of these \text{physically desirable properties.}

\section{A multiscale patch scheme for functional graded beams}
\label{Sec:Patchscheme}

Based on the earlier patch methodology for heterogeneous beams in 2D \citep{Tran2024}, this section develops the application of the method to empower efficient and accurate macroscale predictions of functionally graded viscoelastic beams in~3D.

\subsection{Continuum equations of motion} 
\label{sec:MotionE}

Upon deformation, the elastic stresses arising in a beam are
\begin{equation}
\sigma := 2\mu\varepsilon +\lambda \operatorname{Tr}(\varepsilon )\,I,\label{eq:EStress}
\end{equation}
in which~$\lambda$ and~$\mu$ are Lam\'e constants that characterise the elasticity of the material, and the strain tensor is 
\begin{equation}
\varepsilon(\xv,t):= \frac{1}{2}\left[\nabla \uv +(\nabla \uv )^T\right],\label{eq:EStrain}
\end{equation}
in which $\uv(\xv,t) :=(u,v,w)$ is the displacement vector in directions~\(x,y,z\).

The divergence of the stresses in the material causes unbalanced forces that accelerate the beam. 
On using the Newton's second law, the acceleration of the displacements of a material point is governed by
\begin{equation}
    \rho\DD t\uv=\boldsymbol{\nabla} \cdot\sigma +\fv  \,,\label{E_Disp}
\end{equation} 
in which $\fv $~is the external force density acting on the beam.
For convenience in analysis, the equation of motion \cref{E_Disp} is non-dimensionalized to yield 
\begin{equation}
    \DD {{t^*}}{\uv^*}=\frac{\sigma_0t_0^2}{L_0^2\rho}\boldsymbol{\nabla}^* \cdot\sigma^* +\frac{f_0t_0^2}{L_0\rho}\fv ^* ,\label{eqn:NonDimDisp}
\end{equation}
with the non-dimensional quantities $\uv^*=\uv t_0/L_0$, $t^*=t/t_0$, $\boldsymbol{\nabla}^*=L_0\boldsymbol{\nabla}$, $\sigma^*=\sigma/\sigma_0$ and $\fv ^*=\fv /f_0$.
Then, on choosing characteristic quantities, $L_0=L$, $t_0=L_0/\sqrt{E_{r}/\rho_{r}}$, $\sigma_0=L_0^2\rho_r/t_0^2$, $f_0=L_0\rho_r/t_0^2$, in which $L$~is the length of the beam, and  $E_r$ and~$\rho_r$ are reference values of Young's modulus and the  density of the beam, respectively, the above non-dimensional \pde\ simplifies to
\begin{equation}
    \DD {{t^*}}{\uv^*}=\frac{1}{\rho^*}\left(\boldsymbol{\nabla}^* \cdot\sigma^* +{\fv}^*\right) \,,\label{Eqn:StressDiv}
\end{equation}
in which $\rho^*=\rho/\rho_r$ is dimensionless density.
In the case of functional graded beams, a typical Young's modulus and a typical density are chosen as reference values in the non-dimensionalization.
 

The coordinates system is chosen as follows, and with reference to the undeformed beam: the $x$-axis along the beam longitudinal axis, the $y$-axis and $z$-axis across the beam width and beam thickness, and the origin at the centre of the left-ended cross section, as in the simulation of \cref{FAxiallyGradedBeamSim1}. 
The undeformed beam thickness is~$T$, the beam width is~$W$ and the beam length is~\(L\). 
For the small deformations of linear elasticity, the stress-free boundary conditions
\begin{subequations}\label{eqn:BCs}%
\begin{equation}
    \sigma \cdot\nv_z =0\,,
    \quad\text{and}\quad 
    \sigma \cdot\nv_y =0\,,
\end{equation}
are applied on the top\slash bottom ($z=\pm T/2$) and front\slash back ($y=\pm W/2$) beam surfaces, respectively. 
The left end ($x=0$) of the beam is fixed or clamped, whereas the right end ($x\approx L$) moves freely. 
Accordingly, zero motion applies at the left end,
\begin{equation}
    u=v=w=0 \quad\text{at } x=0\,,
\end{equation}
and a stress-free boundary condition applies over the cross-sectional surface at the free right end,
\begin{equation}
    \sigma \cdot\nv_x =0 \quad\text{at } x=L\,.
\end{equation}
\end{subequations}

\subsection{Microscale spatial discretisation} 
\label{sec:discretisation}

\begin{figure}
\centering
\caption{\label{Fig: MicroGrid}%
Microscale staggered grid in patches in perspective and in three cross-sections.
Displacements evaluated at the middle of the grid edges: $u$,~blue arrow; $v$,~red arrow; $w$,~green arrow. 
Normal stresses evaluated at the grid corners~$\otimes$: $\sigma^{xx}$, $\sigma^{yy}$, $\sigma^{zz}$.  
Shear stresses evaluated at the centre of the grid surfaces: {$\bigcirc$},~$\sigma^{xy}$; $\ominus$,~$\sigma^{xz}$; $\oplus$,~$\sigma^{yz}$. Symbols $\color{teal}\square$: ghost values extrapolated using boundary conditions; $\color{orange}\square$: stress values from stress free boundary conditions; $\color{black}\square$: patch-edge values from patch interpolation.}
\begin{tabular}{@{}l@{}l@{}}
(a) 3D of one micro-grid cube & (b) $xy$-plane \\
\raisebox{6ex}{\input{Figs/3DCell.tex}} & \input{Figs/xyplane.tex} \\
(c) $xz$-plane & (d) $yz$-plane\\
\input{Figs/xzplane.tex} & \input{Figs/yzplane.tex}
\end{tabular}
\end{figure}
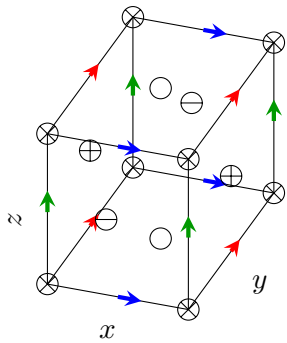

For computational simulation, we must discretise the governing microscale equations~\cref{eq:EStress,eq:EStrain,Eqn:StressDiv,eqn:BCs} in space.
We invoke a microscale staggered discretisation, as illustrated by \cref{Fig: MicroGrid} and akin to several established schemes \citep[e.g.,][]{Virieux84, Virieux86}.

In the Equation-Free scheme, the microscale equations governing the microscale fields are \emph{not} modified in any way: the scheme uses whatever trusted microscale computation a user codes for numerical simulation.
Thus the following discretisation is a standard discretisation and not specific to our multiscale modelling---any valid microscale discretisation scheme is suitable, so we simply chose a common robust microscale discretisation.
\cref{sec:coupling} describes how we compute the microscale in only small \text{sparse patches of the beam.}

We choose to spatially discretise on the microscale staggered grid of \cref{Fig: MicroGrid}, with the grid-points referenced by indices~\(i,j,k\) and spaced by \(\dx,\dy,\dy\) in the \(x,y,z\)~directions, respectively.
Microscale field values, including the displacements and the stresses, are evaluated either at integer microgrid nodes or at half-integer microgrid nodes.
Components~$u,v,w$ of the displacement vector are evaluated at half-integer nodes of the $x,y,z$ directions, respectively, but along integer nodes of the two other directions.  
Normal strains and normal stresses, such as $\varepsilon^{xx}$ and $\sigma^{xx}$, are evaluated at integer nodes for all directions.
Shear strains and shear stresses, such as $\varepsilon^{xy}$ and $\sigma^{xy}$, are evaluated at half-integer nodes in the two subscript directions and integer nodes of the third direction which is perpendicular to the strain or stress.
The temporal derivatives of displacement fields~\cref{Eqn:StressDiv}  are approximated using centred finite differences of stress fields.
For conciseness we define index superscripts~\(\pm\) to denote~\(\pm1/2\) as in \(i\phalf=i+1/2\).
The centred difference operators are~\(\delta_i,\delta_j,\delta_k\), defined such that, for example, \(\delta_i u_{i,j,k}:=u_{i\phalf,j,k}-u_{i\mhalf,j,k}\)\,.
Then, for appropriate ranges of~\(i,j,k\), the spatially discretised equations are 
\begin{subequations}\label{eqs:pdeDisc}%
\begin{align}
     \ddot{u}_{i\phalf,j,k}  & 
     = \frac1{\rho_{i\phalf,j,k}}\left( 
      \frac{\delta_i\sigma^{xx}_{i\phalf,j,k}}{\dx }
     +\frac{\delta_j\sigma^{xy}_{i\phalf,j,k}}{\dy }
     +\frac{\delta_k\sigma^{xz}_{i\phalf,j,k}}{\dz}
     +f^x_{i\phalf,j,k} \right); 
     \label{eq:pdeDiscu}
     \\
    \ddot{v}_{i,j\phalf,k}  & 
    =  \frac1{\rho_{i,j\phalf,k}}\left( 
     \frac{\delta_i\sigma^{xy}_{i,j\phalf,k}}{\dx }
    +\frac{\delta_j\sigma^{yy}_{i,j\phalf,k}}{\dy }
    +\frac{\delta_k\sigma^{yz}_{i,j\phalf,k}}{\dz }
    +f^y_{i,j\phalf,k} \right);
    \\ 
    \ddot{w}_{i,j,k\phalf}  & 
    = \frac1{\rho_{i,j,k\phalf}}\left( 
     \frac{\delta_i\sigma^{xz}_{i,j,k\phalf}}{\dx }
    +\frac{\delta_j\sigma^{yz}_{i,j,k\phalf}}{\dy }
    +\frac{\delta_k\sigma^{zz}_{i,j,k+1}}{\dz }
    + f^z_{i,j,k\phalf} \right).
\end{align}
\end{subequations}
Centred finite differences of the displacement fields approximate the stress field~\cref{eq:EStress} via microscale heterogeneous Lam\'e constants:
\begin{subequations}\label{eqs:StrssEva}%
\begin{align}
    \sigma^{xx}_{i,j,k} & 
    =  (\lambda_{i,j,k}+2\mu_{i,j,k})
    \frac{\delta_iu_{i,j,k}}{\dx }
    +\lambda_{i,j,k}\left(
    \frac{\delta_jv_{i,j,k}}{\dy}
    +\frac{\delta_kw_{i,j,k}}{\dz}
    \right);  
    \\
    \sigma^{yy}_{i,j,k} & 
    =  (\lambda_{i,j,k}+2\mu_{i,j,k})
    \frac{\delta_jv_{i,j,k}}{\dy }
    +\lambda_{i,j,k}\left(
    \frac{\delta_iu_{i,j,k}}{\dx}
    +\frac{\delta_zw_{i,j,k}}{\dz}
    \right);
    \\
    \sigma^{zz}_{i,j,k} & 
    =  (\lambda_{i,j,k}+2\mu_{i,j,k})
    \frac{\delta_kw_{i,j,k}}{\dz }
    +\lambda_{i,j,k}\left(
    \frac{\delta_iu_{i,j,k}}{\dx}
    +\frac{\delta_jv_{i,j,k}}{\dy}
    \right);
    \\
    \sigma^{xy}_{i\phalf,j\phalf,k} & 
    =  \lambda_{i\phalf,j\phalf,k}\left(
    \frac{\delta_ju_{i\phalf,j\phalf,k}}{\dy }
    +\frac{\delta_iv_{i\phalf,j\phalf,k}}{\dx } \right);
    \\
    \sigma^{yz}_{i,j\phalf,k\phalf} & 
    =  \lambda_{i,j\phalf,k\phalf}\left( 
    \frac{\delta_jw_{i,j\phalf,k\phalf}}{\dy }
    +\frac{\delta_kv_{i,j\phalf,k\phalf}}{\dz } \right);
    \\
    \sigma^{xz}_{i\phalf,j,k\phalf} & 
    =  \lambda_{i\phalf,j,k\phalf}\left(
    \frac{\delta_ku_{i\phalf,j,k\phalf}}{\dz }
    +\frac{\delta_iw_{i\phalf,j,k\phalf}}{\dx } \right).
\end{align}
\end{subequations}

The physical elastic stress-free boundary conditions~\cref{eqn:BCs} on the top\slash bottom and front\slash back beam surfaces (and as used for the patches in \cref{FAxiallyGradedBeamSim1}) are implemented via the ghost values of~$v$ and~$w$ in \cref{Fig: MicroGrid}(b)--(d).
In the discretisation of the cross-section of the beam, \cref{FAxiallyGradedBeamSim1}(d), we impose a microgrid of \(n_y\times n_z\) points indexed by~\(j,k\).
The microgrid index~\(i\) for the \(x\)-direction is specified in \cref{sec:coupling} where we also discuss  constructing patches along the beam.
\begin{itemize}
\item On the front surface \(y=-W/2\) (index $j=1$), for every appropriate~\(i\): the normal stress $\sigma^{yy}=0$ for interior nodes \(k=2:n_z-1\); and $\sigma^{yy}=\sigma^{zz}=0$ for edge nodes \(k=1,n_z\).  Hence 
\begin{subequations}\label{eqs:bcleft}%
\begin{align*}
     v_{i,1\mhalf,1} &= v_{i,1\phalf,1}
     +\frac{\dy}{\dx}\frac{\lambda_{i,1,1}}{2\left(\mu_{i,1,1}+\lambda_{i,1,1}\right)}\delta_i u_{i,1,1}; \\
     v_{i,1\mhalf,k} &= v_{i,1^+,k}
     +\frac{\dy\lambda_{i,1,k}}{2\mu_{i,1,k}+\lambda_{i,1,k}}
     \left(\frac{\delta_iu_{i,1,k}}{\dx}
     +\frac{\delta_k w_{i,1,k}}{\dz}\right);
     \\
      v_{i,1\mhalf,n_z} &= v_{i,1^+,n_z}
      +\frac{\dy}{\dx}\frac{\lambda_{i,1,n_z}}
      {2(\mu_{i,1,n_z}+\lambda_{i,1,n_z})}\delta_iu_{i,1,n_z}.
\end{align*}
\end{subequations}

\item On the back surface  \(y=+W/2\) (index $j=n_y$), for every appropriate~\(i\): the normal stress $\sigma^{yy}=0$ for interior nodes \(k=2:n_z-1\); and $\sigma^{yy}=\sigma^{zz}=0$ for edge nodes \(k=1,n_z\).  Hence 
\begin{subequations}\label{eqs:bcright}%
\begin{align*}
     v_{i,n_y\phalf,1} &= v_{i,n_y\mhalf,1}-\frac{\dy}{\dx}\frac{\lambda_{i,n_y,1}}{2\left(\mu_{i,n_y,1}+\lambda_{i,n_y,1}\right)}\delta_iu_{i,n_y,1};\\
     v_{i,n_y\phalf,k} &= v_{i,n_y\mhalf,k}-\frac{\dy\lambda_{i,n_y,k}}{2\mu_{i,n_y,k}+\lambda_{i,n_y,k}}\left(\frac{\delta_iu_{i,n_y,k}}{\dx}+\frac{\delta_kw_{i,n_y,k}}{\dz}\right);\nonumber\\
      v_{i,n_y\phalf,n_z} &= v_{i,n_y\mhalf,n_z}-\frac{\dy}{\dx}\frac{\lambda_{i,n_y,n_z}}{2\left(\mu_{i,n_y,n_z}+\lambda_{i,n_y,n_z}\right)}\delta_iu_{i,n_y,n_z}\,.
\end{align*}
\end{subequations}
\item On the bottom surface \(z=-T/2\)  (index $k=1$), for every appropriate~\(i\): the normal stress $\sigma^{zz}=0$ for interior nodes \(j=2:n_y-1\); and $\sigma^{yy}=\sigma^{zz}=0$ for edge nodes \(j=1,n_y\).  Hence 
\begin{subequations}\label{eqs:bcbot}%
\begin{align*}
     w_{i,1,1\mhalf} &= w_{i,1,1\phalf}+\frac{\dz}{\dx}\frac{\lambda_{i,1,1}}{2\left(\mu_{i,1,1}+\lambda_{i,1,1}\right)}\delta_iu_{i,1,1};\\
     w_{i,j,1\mhalf} &= w_{i,j,1\phalf}+\frac{\dz\lambda_{i,j,1}}{2\mu_{i,j,1}+\lambda_{i,j,1}}\left(\frac{\delta_iu_{i,j,1}}{\dx}+\frac{\delta_jv_{i,j,1}}{\dy}\right);\nonumber\\
      w_{i,n_y,1\mhalf} &= w_{i,n_y,1\phalf}+\frac{\dz}{\dx}\frac{\lambda_{i,n_y,1}}{2\left(\mu_{i,n_y,1}+\lambda_{i,n_y,1}\right)}\delta_iu_{i,n_y,1}\,.
\end{align*}
\end{subequations}

\item On the top surface \(z=+T/2\)  (index $k=n_z$), for every appropriate~\(i\): the normal stress $\sigma^{zz}=0$ for interior nodes \(j=2:n_y-1\); and $\sigma^{yy}=\sigma^{zz}=0$ for edge nodes \(j=1,n_y\). Hence 
\begin{subequations}\label{eqs:bctop}%
\begin{align*}
     w_{i,1,n_z\phalf} &= w_{i,1,n_z\mhalf}-\frac{\dz}{\dx}\frac{\lambda_{i,1,n_z}}{2\left(\mu_{i,1,n_z}+\lambda_{i,1,n_z}\right)}\delta_iu_{i,1,n_z};\\
     w_{i,j,n_z\phalf} &= w_{i,j,n_z\mhalf}-\frac{\dz\lambda_{i,j,n_z}}{2\mu_{i,j,n_z}+\lambda_{i,j,n_z}}\left(\frac{\delta_iu_{i,j,n_z}}{\dx}+\frac{\delta_jv_{i,j,n_z}}{\dy}\right);\nonumber\\
      w_{i,n_y,n_z\phalf} &= w_{i,n_y,n_z\mhalf}-\frac{\dz}{\dx}\frac{\lambda_{i,n_y,n_z}}{2\left(\mu_{i,n_y,n_z}+\lambda_{i,n_y,n_z}\right)}\delta_iu_{i,n_y,n_z}\,.
\end{align*}

\item The zero shear stress condition on the top\slash bottom\slash front\slash back beam surfaces applies directly to the grid nodes over the surfaces, for \(j=2:n_y-1\)\,, \(k=2:n_z-1\)\,, and for every appropriate~\(i\):
    \begin{align*}
        \sigma^{xy}_{i\phalf,1\mhalf,k}&=-\sigma^{xy}_{i\phalf,1\phalf,k}\,;&
        \sigma^{xy}_{i\phalf,n_y\phalf,k}&=-\sigma^{xy}_{i\phalf,n_y\mhalf,k}\,;\\
        \sigma^{xz}_{i\phalf,j,1\mhalf}&=-\sigma^{xz}_{i\phalf,j,1\phalf}\,;&
        \sigma^{xz}_{i\phalf,j,n_z\phalf}&=-\sigma^{xz}_{i\phalf,j,n_z\mhalf}\,;\\
        \sigma^{yz}_{i,j\phalf,1\mhalf}&=-\sigma^{yz}_{i,j\phalf,1\phalf}\,;&
        \sigma^{yz}_{i,j\phalf,n_z\phalf}&=-\sigma^{yz}_{i,j\phalf,n_z\mhalf}\,;\\
        \sigma^{yz}_{i,1\mhalf,k\phalf}&=-\sigma^{yz}_{i,1\phalf,k\phalf}\,;&
        \sigma^{yz}_{i,n_y\phalf,k\phalf}&=-\sigma^{yz}_{i,n_y\mhalf,k\phalf}\,.
    \end{align*}
\end{subequations}
\end{itemize}
The other physical boundary conditions apply at the two ends of the beam at \(x=0,L\).
The left end of the beam is fixed, so the rigid zero displacement applies to the leftmost microgrid points: 
\begin{align*}
    u_{1\phalf,j,k}^1&=0\quad \text{for } j=1:n_y,\, k=1:n_z;\\
    v_{1,j\phalf,k}^1 &= 0\quad \text{for } j=1:n_y-1,\, k=1:n_z;\\
    w_{1,j,k\phalf}^1 &=0\quad \text{for } j=1:n_y,\, k=1:n_z-1.
\end{align*}    
For the free end on the right, the stress-free condition applies directly to the rightmost microgrid points: 
\begin{align*}
    \sigma_{n_x,j,k}^{xx,N} 
    &=-\sigma_{n_x-1,j,k}^{xx,N}\quad \text{for } j=1:n_y,\, k=1:n_z;\\
    \sigma_{n_x\mhalf,j\phalf,k}^{xy,N} &=0 \quad \text{for } j=1:n_y-1,\, k=1,n_z;\\
    \sigma_{n_x\mhalf,j,k\phalf}^{xz,N} &=0\quad \text{for } j=1:n_y,\, k=1:n_z-1.
\end{align*}

\paragraph{Viscoelastic dissipation}
Physically, weak dissipation via boundary friction with the ambient environment, and\slash or sound energy radiated to the surrounds, and/or internal viscoelasticity would typically occur. 
To include within the scope of this exploration some representative physical effects of such dissipation, we include a simple Kelvin--Voigt viscoelastic term,~$+\eta \nabla^2\dot{\uv}$, into the non-dimensional \pde~\cref{Eqn:StressDiv}.
This viscoelasticity adds into the microscale  discretisation~\cref{eq:pdeDiscu}, for~\(\ddot u_{i\phalf,j,k}\),  
\begin{align}
\cdots&+\eta\left(
\frac{\delta_i^2u_{i\phalf,j,k}}{\dx ^2}
+\frac{\delta_j^2u_{i\phalf,j,k}}{\dy ^2}
+\frac{\delta_k^2u_{i\phalf,j,k}}{\dz ^2}
\right),
\label{eqn:viscous}
\end{align}
in which $\delta_i^2u_{i\phalf,j,k}:=\delta_i\left(\delta_iu_{i\phalf,j,k}\right)=\delta_iu_{i+1,j,k}-\delta_iu_{i,j,k}$. 
Discretization of the dissipation terms for \(\ddot v_{i,j\phalf,k}\) and \(\ddot w_{i,j,k\phalf}\) is calculated similarly.
This Kelvin--Voigt viscoelasticity is also representative of a range of other phenomenological dissipation with non-dimensional strength~\(\eta\).
Typically we choose small viscoelastic effects, \(\eta=10^{-3}\). 

\subsection{Resolve microscale only on small sparse patches in space}
\label{sec:coupling}

We compute the discretised viscoelastic equations~\cref{eqs:pdeDisc}--\cref{eqn:viscous} of 3D beams within $N$~well-separated, small equi-sized patches, each of length \(h\ll L\), width~\(W\) and thickness~$T$.
\cref{FAxiallyGradedBeamSim1} shows an example for a beam of (non-dimensional) length \(L=1\) and thicknesses \(W=T=0.09\) computed only within $N=9$ patches each of \(x\)-length \(h=0.036\)\,.
These patches correspond to Representative Volume Elements in other multiscale methods, but importantly we obtain very controllable high accuracy with no subjective averaging, no oversampling regions, no buffer regions, and no guessed fast/slow variables.

The macroscale spacing~$H$ between the patches is chosen to resolve the required macroscale dynamics of the beams.
Microscale computations are performed only within the patches.
In each patch, visible in the grids of \cref{FAxiallyGradedBeamSim1,Fig: MicroGrid}, let the number of microgrid points within each patch in the $x$-dimension be denoted by~$n_x$ (recall the beam's \(yz\)-cross-sections have a $n_y\times n_z$ microgrid). 
When unstressed, the microgrid nodes are equi-spaced with $\dx:=h/(n_x-1)$, $\dy:=W/(n_y-1)$ and $\dz:=T/(n_z-1)$ in the~$x$, $y$ and~$z$~dimensions.  
For example, \cref{Fig: MicroGrid}(b)--(d) shows sections through one patch with a \(n_x\times n_y\times n_z =7\times 6\times 6\) microgrid.     

\begin{figure}
\centering
\caption{\label{figpatchscheme}%
Spatial patches, indexed by~\(I\), with a macroscale spacing~\(H\gg h\). 
Each patch interior (green) is bounded by patch-edge values (black boxes).
Interpolating left\slash right next-to-edge values (blue/red arrows) from neighbouring patches determines the edge values on the right\slash left of each patch~\(I\) \citep{Bunder2020a}.   
Interpolating from patches \(I-2,\ldots,I+2\) to the edges of the $I$th~patch (as shown) determines the macroscale to errors~\Ord{H^4}. 
}
\input{Figs/figpatchscheme.tex} 
\end{figure}
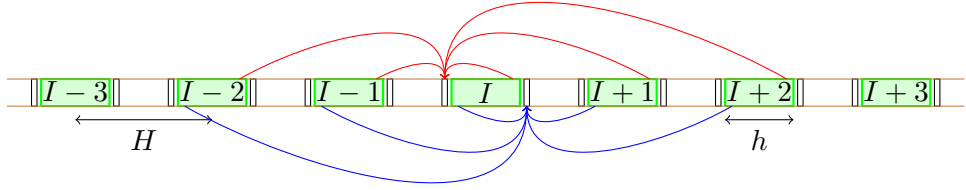%
\paragraph{Patch coupling for macroscale dynamics}
The patches are coupled across unsimulated space to capture the macroscale dynamics.
Although there are various good methods to couple the patches \citep[e.g.,][]{Roberts2023a}, we prioritize patch-coupling methods that
preserve self-adjoint symmetries of the microscale model \citep{Bunder2020a}, and hence better preserve numerical stability in dynamic simulations.   
\cref{figpatchscheme} schematically represents the self-adjoint coupling method.
Here, the interpolated quantities are the displacement fields~\((u,v,w)\) at the left and right faces of each patch.
For each patch index~\(I\) and cross-beam micro-grid indices~\(j\) and~$k$, displacement values on the right face of the patch, \(u_{n_x\phalf,j,k}^I\), \(v_{n_x,j\phalf,k}^I\), and~\(w_{n_x,j,k\phalf}^I\) ($\square$\,s in \cref{Fig: MicroGrid}), are determined from interpolation of the displacement values, \(u_{2\phalf,j,k}^J\), \(v_{2,j\phalf,k}^J\) and~\(w_{2,j,k\phalf}^J\) respectively, of the left next-to-face grid-points of neighbouring patches indexed by~\(J\). 
Correspondingly, displacement fields on the left face of a patch \(u_{1\phalf,j,k}^I\), \(v_{1,j\phalf,k}^I\) and~\(w_{1,j,k\phalf}^I\) are determined from interpolation of displacement fields \(u_{n_x\mhalf,j,k}^J\), \(v_{n_x-1,j\phalf,k}^J\) and~\(w_{n_x-1,j,k\phalf}^J\), respectively, of the right next-to-face grid-points of neighbouring patches~\(J\). 
These interpolated patch face values then affect predictions of the patch scheme via the stress field equations~\cref{eqs:StrssEva}. 
Here we apply a polynomial interpolation using the next-to-face displacement values of the patches.
An order~\(P\) polynomial interpolation means that the left\slash right face field of the \(I\)th~patch is determined from a polynomial interpolation of right\slash left  next-to-edge fields of patches indexed by \(J=I-P/2,\ldots, I+P/2\).

\paragraph{Patch size}
An advantage of the patch scheme is that microscale computations are only carried out within patches that cover a small portion of the full-domain. 
A crucial decision in applying the patch scheme is the size of the patches in space: the smaller the patches the greater the computational savings, but not all small patches will provide accurate simulations.
The \emph{patch ratio} defined by \(r:=h/H<1\) characterises the small fraction of the beam within which the microscale computations are performed.
Herein, we assume slender beams so choose each patch to extend across the full cross-sectional area of the beam and hence all patches are arranged along the longitudinal axis of the beam.
The smallest patch size that also provides accurate simulations is when \(h\)~corresponds to one microscale period \citep{Bunder2017}, and this is the patch size chosen herein.

In the case of functionally graded, near-periodic, microscale heterogeneity \cite[e.g.,][]{Anthoine2010, Shahraki2020, Shahbaziana2022}, the on-the-fly computations of the patch scheme automatically adapt to any functional graduations that occur on the macroscale.

\section{Cross-sectional graded beam}
\label{Sec:CrSecGB}


This section develops the application of the patch scheme to 3D heterogeneous beams in order to study static equilibrium and dynamics of cross-sectional graded cantilever beams of~\(\Al/\Si\C\).
Such beams were experimentally studied by \cite{Kapu2008}, hereafter denoted by \kbk. 
The numerical results on maximum deflection and natural frequencies of the beams are here compared with the experimental data reported by \kbk.

Specifically, this section reports on two beam structures.
\begin{itemize}
  \item A three-layer~\Al/\Si\C\ beam with each layer of the beam made from a mixture of~\Al\ and~\Si\C\ with volume fraction ratios $5:0$~in the bottom layer, $4:1$~in the middle, and $3:2$~in the top.  
  \item A five-layer~\Al/\Si\C\ beam with volume fraction ratios of~\Al\ to~\Si\C\  from the bottom to the top layers of~\(10:0\), \(9:1\), \(8:2\), \(7:3\), and~\(6:4\), respectively.
\end{itemize}
For the \Al/\Si\C~composite with volume fractions~\(V_m\) and~\(V_c\), we adopt the following mixing rule for Poison's ratio and the density of the mixing layers (\kbk):
\begin{equation}
      \rho_{l}=\rho_mV_m+\rho_cV_c\,,\qquad
      \nu_{l}=\nu_mV_m+\nu_cV_c\,,
      \label{eq:LinMixRule}
\end{equation}
in which subscript~\(l\) denotes a quantity of the \Al/\Si\C~composite, and subscripts~\(m\) and~\(c\) denote quantities of metal~\Al\ and ceramic~\Si\C\ components, respectively.  
For the Young's modulus of the mixture, we apply a semi-empirical mixing rule (\kbk):  
\begin{equation}
   E=\frac{V_mE_m\frac{q+E_c}{q+E_m}+V_cE_c}{V_m\frac{q+E_c}{q+E_m}+V_c}\,,\label{eq:EMixRule}
\end{equation}
in which experiments (\kbk) determine the parameter $q=91.6\opn{GPa}$ for~\Al/\Si\C\ from the stress-strain relation of the two mixing components.
\cref{Tab:BeamDatCSGB} summarises the physical properties of the materials. 

\begin{table}
\caption{\label{Tab:BeamDatCSGB}Properties of the materials and the cantilever beams in the experiments of \protect\cite{Kapu2008} (\kbk).}
\centering
\begin{tabular}{@{}l@{}l@{}}
\hline  
    \(\Al\)&\\
\hline
     Young's modulus $E_{\Al}$  & $67\opn{GPa}$  \\
     Density  $\rho_{\Al}$      & $2700\opn{kg}\opn{m}^{-3}$\\
     Poisson ratio $\nu_{\Al}$  & $0.33$ \\   
\hline           
     \(\Si\C\)&\\
\hline
     Young's modulus $E_{\Si\C}$ & $302\opn{GPa}$  \\
     Density $\rho_{\Si\C}$      & $3200\opn{kg}\opn{m}^{-3}$\\
     Poisson ratio $\nu_{\Si\C}$ & $0.17$ \\
\hline                       
    three-layer \(\Al/\Si\C\) beam &\\
\hline
    Length $L$ & $105\opn{mm}$\\
    Width $W$                & $15\opn{mm}$\\
    Layer thickness $h_l$    & $3\opn{mm}$\\  
    Number of layers $n_l$   & $3$\\
    Volume ratio ($V_m:V_c$) & $100:0$, $80:20$, $60:40$\\ 
\hline 
    five-layer \(\Al/\Si\C\) beam& \\
\hline
    Length $L$ & $110\opn{mm}$\\
    Width $W$              & $15\opn{mm}$\\
    Layer thickness $h_l$  & $2\opn{mm}$\\  
    Number of layers $n_l$ & $5$\\
    Volume ratio ($V_m:V_c$)& $100:0$, $90:10$, $80:20$, $70:30$, $60:40$ \\
\hline 
\end{tabular}
\end{table}

\cref{Tab:BeamDatCSGB} also summarises the beam dimensions.
The beams were manufactured to be~$125\opn{mm}$ long. 
However, in the experiments of \kbk\ the distance from the clamped\slash fixed point to the tip of the free end is smaller than the manufactured length: $105\opn{mm}$ and~$110\opn{mm}$ for the three-layer and five-layer \Al/\Si\C\ beams, respectively. 
These distance are taken to be the length~\(L\) of the beam for our computations.
The beam thickness depends upon layer thickness and the number of layers.
The layer thickness of the three-layer beam and the five-layer beam is~$3\opn{mm}$ and~$2\opn{mm}$, respectively.
Hence, the thicknesses of the three-layer and five-layer beams are~$9\opn{mm}$ and~$10\opn{mm}$, respectively.
In the static equilibrium study, a static load $f=147\opn{N}$ is applied at the free end of the beams. 

We implement patch dynamics for~$9$ and~$17$~patches along the beam.  
Each of the patches is defined by a microgrid with $(n_x,n_y)= (7,5)$ and \(n_z=6\) and \(10\) for the three-layered and five-layered beams, respectively.
Accordingly, the microscale spatial resolutions $(\dx,\dy,\dz)=(0.0061,0.0341,0.0143)$ for the three-layer beam and $(0.0061,0.0341,0.0091)$ for the five-layer beam.
With this fixed spatial resolution, $9$~patches and $17$~patches cover~$27\%$ and~$52\%$ of the beam domain, respectively.
Experimental data from \kbk\ and from full-domain computations in \cref{Tab:MaxDefCSGB,Tab:FreCSGB} validate the predictions of the patch scheme.

\begin{figure}\centering
\caption{\label{fig:ErrorAlSiC} Axial displacement of the three-layer and five-layer \Al/\Si\C\ beams computed using $9$~patches, $17$~patches, and the full beam domain. 
The errors are relative to the maximum displacement in the corresponding full-domain computations.}
\begin{tabular}{@{}l@{}l@{}}
(a) 3-layered beam&\\
\inPlot{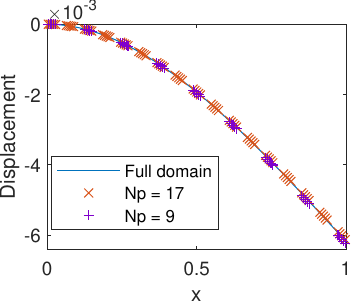}&
\inPlot{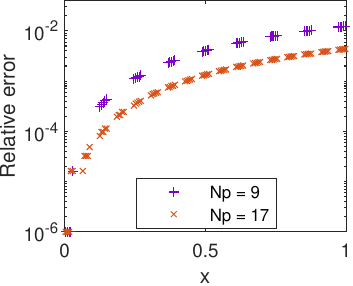}\\
(b) 5-layered beam &\\
\inPlot{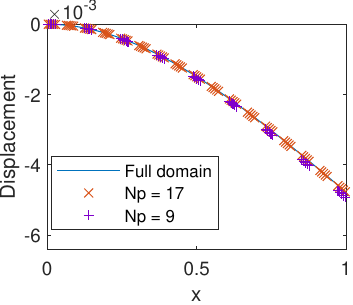}&
\inPlot{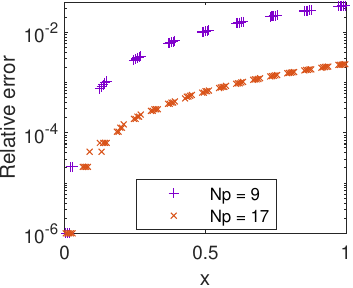}
\end{tabular}
\end{figure}

\cref{fig:ErrorAlSiC} displays computed beam displacements along the beam axis for the equilibrium displacements of the three-layer and five-layer \Al/\Si\C\ beams.
The displacements are measured along the center line of the beams at interior microgrid nodes of the patches. 
The blue solid curve represents the axial displacements of the beam simulated in the full-domain, whereas `$\times$'~and~`$+$' are for computations of~$17$ and~$9$~patches, respectively.
The displacement values in the $17$~patches and $9$~patches cases align very well with the full-domain curves, and the deflections in the patches give a clear indication of the full beam deflection in this static equilibrium state. 
\cref{fig:ErrorAlSiC} also shows errors relative to the maximum displacement at the beam tip of the full-domain computation.
The errors consistently increase along the beam axis from the fixed end to the free end, and are larger for $9$~patches than for $17$~patches. 
The maximum error in the three-layer beam is~$1.1\%$ for $9$~patches  and~$0.1\%$ for $17$~patches, whereas in the five-layer beam  they are~$1.9\%$ and~$0.2\%$, respectively.

\cref{Tab:MaxDefCSGB} summarizes the maximum deflection values predicted by the computations together with the corresponding experimental values by \kbk.
\kbk\ repeated each experiment three times, so the table provides their observed range of maximum deflection, and they are non-dimensionalized by the beam length (\cref{Tab:BeamDatCSGB}). 
\cref{Tab:MaxDefCSGB} gives two important conclusions. 
Firstly, there is good agreement between the computational values and the experimental values for both beam cases and for both patch numbers. 
Secondly, the errors increase by about~$1\%$ as the patch number decreases from~$17$ to~$9$. 
These results are impressive as the microscale computations are carried out in only a half and a quarter of the beam domain.
The agreement between the computational and experimental values for 3D functionally graded beams indicates that the computational patch scheme is effective in evaluating static deflection of cross-sectional graded beams.
\begin{table}
\caption{\label{Tab:MaxDefCSGB}Maximum deflection of the \(\Al/\Si\C\) cantilever beams obtained from the patch computations and measured in the \kbk\ experiments.}
\centering
\(\begin{array}{p{6em}cc}
\hline    
    & \text{three-layer }\Al/\Si\C & \text{five-layer }\Al/\Si\C \\   
\hline    
   $9$~patches   & 6.23\xten{-3} & 4.85\xten{-3}\\           
   $17$~patches  & 6.13\xten{-3} & 4.77\xten{-3}\\         
   full-domain  & 6.16\xten{-3} & 4.76\xten{-3}  \\
   Experiment  & 6.2-7.4\xten{-3} & 5.5-6.3\xten{-3} \\
\hline    
\end{array}\)
\end{table}

The computational efficiency of the patch scheme is assessed via the computation time. 
The equilibria were all computed on a desktop computer using \matlab's \verb|fsolve()| due to its robustness and simple coding, albeit at the cost of generally longer compute times than is likely when using conjugate gradients and incomplete \textsc{lu}-preconditioning.
For the three-layer beam, the computation time for the full-domain computation is~$95\opn{s}$, whereas it is~$32\opn{s}$ and~$12\opn{s}$ for  $17$~patches and $9$~patches, respectively.
For the five-layer beam the compute times are~$357\opn{s}$, $69\opn{s}$, and~$21\opn{s}$, respectively. 
\cref{Tab:CompTimeCSGB} summarizes the static equilibrium computation times. 
For the three-layer beam, using $17$~patches and $9$~patches reduces the full-domain computation time by factors of three and eight, respectively.
For the five-layer beam, the patch scheme reduces the time by factors of five and seventeen, respectively.
This demonstrates that the patch scheme is significantly faster than a full-domain simulation.
Much larger gains in computational speed are obtained in other patch dynamics applications, for example \cite{Divahar2022b} obtain speed-ups of up to 100,000 in 2D wave simulations. 
\begin{table}
\caption{\label{Tab:CompTimeCSGB}Computational time in seconds for beam deflection and eigenvalues of the Jacobian matrix.   
Here the patch scheme takes as little as~6\% of the time of a full-domain computation. }
\centering
\begin{tabular}{|l|ccc|}
\hline
   & $9$~patches & $17$~patches & full-domain \\
\hline
\hline
  & \multicolumn{3}{c|}{Deflection}\\
\hline  
 three-layer \Al/\Si\C\ & $12$ & $32$ & $95$\\
 five-layer \Al/\Si\C\ & $21$ & $69$ & $357$\\
\hline
\hline
  & \multicolumn{3}{c|}{Eigenvalues and eigenvectors}\\
\hline
 three-layer \Al/\Si\C\ & $5$ & $17$ & $46$ \\
 five-layer \Al/\Si\C\ & $10$ & $32$ & $143$  \\
\hline
\end{tabular}
\end{table}

\begin{figure}\centering
\caption{\label{fig:EValue} Leading 4\% of the complex-valued eigenvalues of the Jacobian matrices for the three-layer and five-layer beams 
from the $9$~patch and $17$~patch computations. 
The~$\circ$, {\tiny$\square$}, $\triangleright$, and~{\small$\vartriangle$} symbols indicate, respectively, the $z$-bending modes,  $y$-bending modes, torsion modes, and axial compression modes: red, green, and yellow denote the first, second and third modes.  The axes are quasi-logarithmic.}
\begin{tabular}{@{}l@{}l@{}}
(a) three layers and $9$~patches & 
(b) three layers and $17$~patches\\
\inPlot{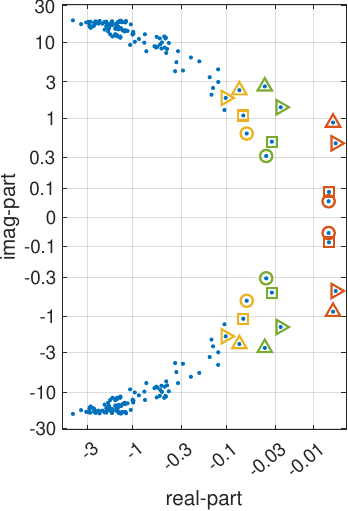} & 
\inPlot{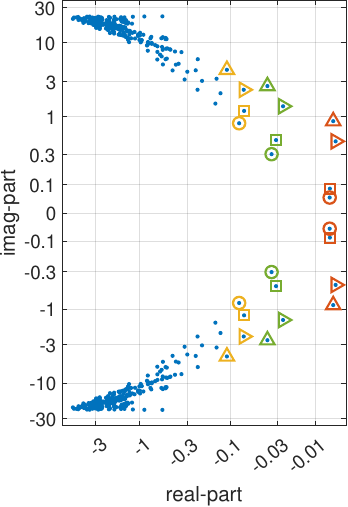}    \\
(c) five layers and $9$~patches & 
(d) five layers and $17$~patches\\
\inPlot{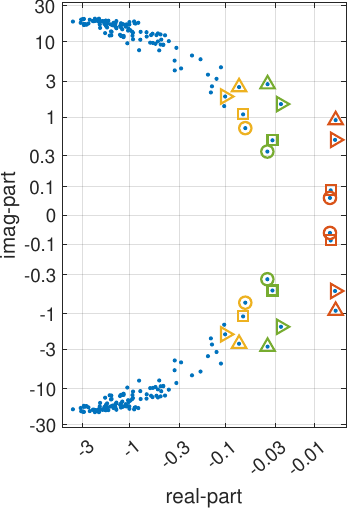} & 
\inPlot{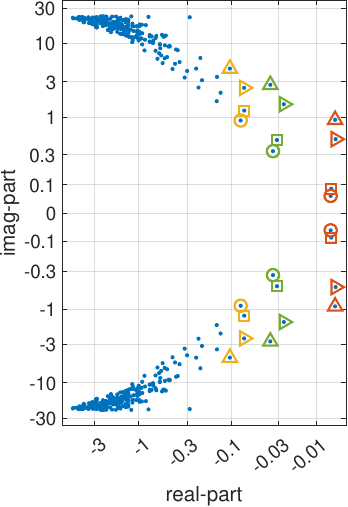}\\
\end{tabular}
\end{figure}

\begin{figure}\centering
\caption{\label{fig:EModes5AlsiC917} The first three \(z\)-bending modes of the five-layer \Al/\Si\C\ beam simulated using $9$~patches (left) and $17$~patches (right). 
Mesh colour is the axially normal stress~$\sigma^{xx}$.}
\begin{tabular}{@{}l@{}l@{}}
\raisebox{\height}{\parbox{0.5\textwidth}{%
(a) The first $z$-bending mode\hfil\\
\inPlot{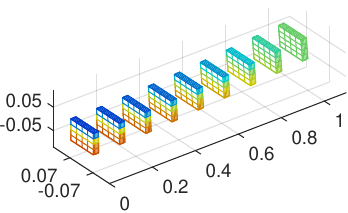}}} & 
\inPlot{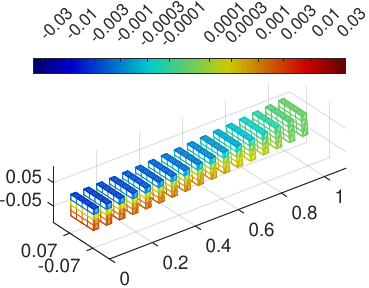}\\ 
\raisebox{\height}{\parbox{0.5\textwidth}{%
(b) The second $z$-bending mode\hfil\\
\inPlot{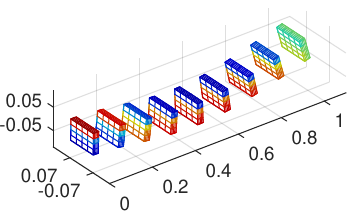}}} & 
\inPlot{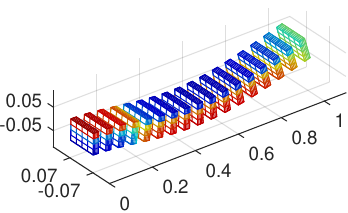}\\
\raisebox{\height}{\parbox{0.5\textwidth}{%
(c) The third $z$-bending mode\hfil\\
\inPlot{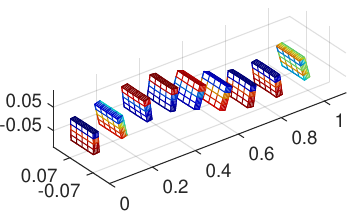}}} & 
\inPlot{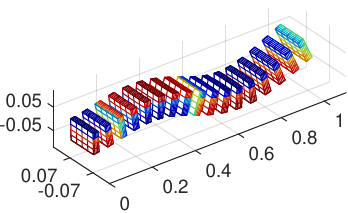}\\
\end{tabular} 
\end{figure}

\begin{figure}\centering
\caption{\label{fig:OEModes5AlsiC917} The first two natural vibration modes of the five-layer \Al/\Si\C\ beam simulated using $17$~patches: $y$-bending modes (1st row), torsion modes (second row), and compression modes (third row).  
The left and right columns correspond to the first and second modes, respectively.
Mesh colour is the axially normal stress~$\sigma^{xx}$.}
\begin{tabular}{@{}l@{}l@{}}
\raisebox{\height}{\parbox{0.5\textwidth}{(a) The first two $y$-bending modes\hfil\\
\inPlot{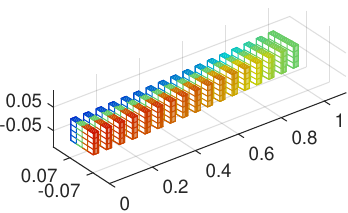}}} & \inPlot{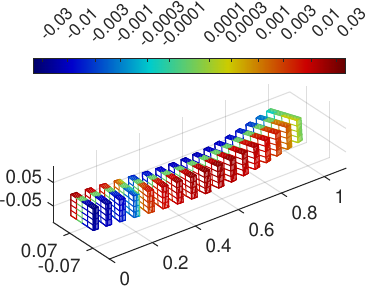}\\
(b) The first two torsion modes\hfil \\ 
\inPlot{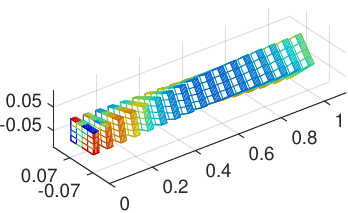} & \inPlot{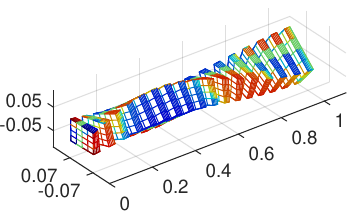}\\
(c) The first two compression modes\hfil\\
\inPlot{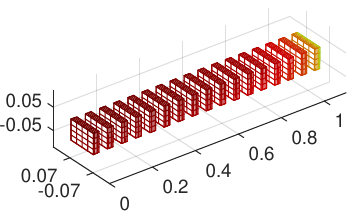} & \inPlot{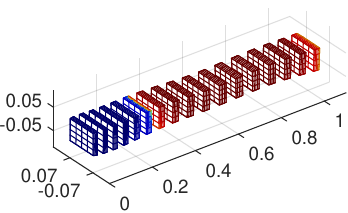}\\
\end{tabular} 
\end{figure}

The above discussion addresses equilibria of the beam.  
However, the dynamics of a beam are also crucial in many physical scenarios, and theory underpinning the patch scheme supports its application to general dynamical systems \cite[e.g.,][]{Roberts06d, Roberts2011a, Bunder2020a}.
Here, the governing viscoelastic equations are linear, so we know that an eigenvalue-eigenvector decomposition provides an understandable general solution to the dynamics of the original equations as well as to the patch versions. 
In each case, eigenvalues and eigenvectors are computed from a numerically computed Jacobian matrix of the equations of motion (\cref{Sec:Patchscheme}), as done by \cite{Tran2024}.
The vibration modes of the beam come from the eigenvectors.
The imaginary part of an eigenvalue is the angular frequency of the corresponding eigenmode, whereas the magnitude of the real part quantifies the damping of the eigenmode due to viscoelastic dissipation.
\cref{fig:EValue} plots the eigenvalues of the beams computed using $9$~patches and $17$~patches.
In these plots, the smallest magnitude real-part eigenvalues are highlighted as they are the slowest to decay so dominate the long-term dynamics.
The eigenvalues are classified into four types of vibration modes: $z$-bending modes, $y$-bending modes, torsion modes in the $yz$-plane, and axial compression modes.
Comparing the eigenvalues of $9$~and~$17$ patches, slight differences are seen in the third modes of the four vibrating types.
\cref{fig:EModes5AlsiC917,fig:OEModes5AlsiC917} illustrate typical examples of the four different modes.  

\cref{fig:EValue} reflects several important physical beam properties. 
Firstly, mode frequencies in ascending order are  $z$-bending,  $y$-bending, torsion, and compression.
Secondly, for each of the four types of vibration modes, the viscoelastic damping increases as the frequency increases. 
Finally, the three-layer and five-layer beams have very similar eigenvalues for the first three modes of the four types of vibration. 
Consequently, vibrations of the two beams are expected to be similar under the same conditions.      
 
\cref{Tab:FreCSGB} lists values of the first three angular $z$-bending frequencies of the beams computed using $9$~patches, $17$~patches, and the full-domain.
Non-dimensional frequencies obtained from experiments by \kbk\ are also included for comparison. 
The relationship between a physical frequency from the experiments and a non-dimensional angular frequency is $\omega^*=2\pi t_0f$, with $t_0=L_0/\sqrt{E_{\max}/\rho_{\max}}$ the characteristic time.  
In the \kbk\ experiments, beam vibrations were generated and  displacements were recorded over a period of time. 
The displacement data in the time domain was transformed to the frequency domain using a Fourier transform to obtain the natural frequencies of the beam vibration. 
\cref{Tab:FreCSGB} shows the ranges of frequencies observed in the \kbk\ experiments.
The computed and experimental frequencies are quite similar.
This verifies that the patch dynamics algorithm for elasticity not only calculates accurate static deformations of beams, but also correctly determines the dynamics of beams.   
\begin{table}
\caption{\label{Tab:FreCSGB}The first three non-dimensional angular frequencies~$\omega^*$ (to three decimal places) of the bending mode of the cross-sectional graded beams simulated using $9$~patches, $17$~patches, and the full-domain,  compared with those obtained from the \kbk\ experiments. }
\begin{equation*}
\begin{array}{|c|cccc|}
\hline
\text{Modes} & \text{$9$~patches} & \text{$17$~patches} & \text{full-domain} & \text{Experiment}\\
\hline
\hline
 & \multicolumn{4}{c|}{\text{three-layer }\Al/\Si\C}\\      
\hline
1 & 0.050 & 0.051 & 0.051 & 0.055-0.060\\
2 & 0.313 & 0.307 & 0.306 & 0.340-0.367\\
3 & 0.875 & 0.873 & 0.873 & 0.922-0.988\\
\hline
\hline
 & \multicolumn{4}{c|}{\text{five-layer }\Al/\Si\C}\\      
\hline
1 & 0.056 & 0.057 & 0.057 & 0.051-0.054\\
2 & 0.345 & 0.340 & 0.339 & 0.312-0.332\\
3 & 0.916 & 0.915 & 0.915 & 0.838-0.883\\
\hline
\end{array}
\end{equation*}
\end{table}

\begin{figure}
\centering
\caption{\label{fig:DefAlB} Displacement measured at the free tip of the three-layer and five-layer \Al/\Si\C\ beams over the non-dimensional simulation time.}
\begin{flushleft} (a) three-layer beam\\ \end{flushleft}
\inPlot{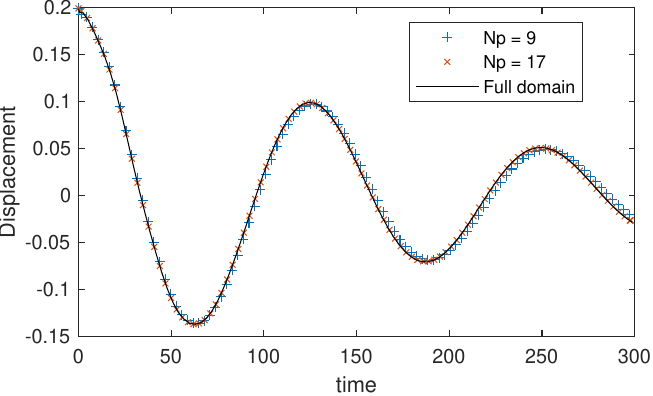} \\
\begin{flushleft} b) five-layer beam\\ \end{flushleft}
\inPlot{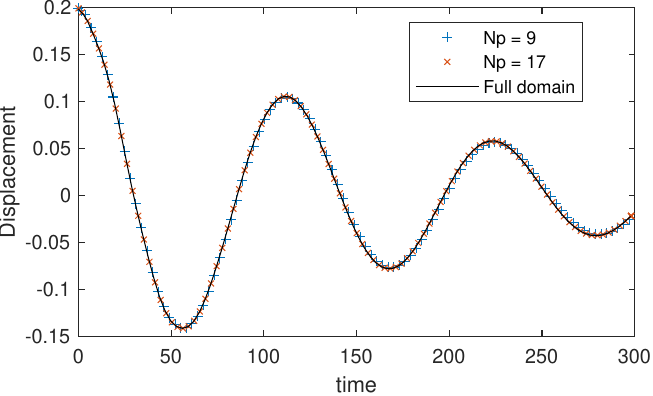}  
\end{figure}

The vibration of the cross-sectional graded~\Al-\Si\C\ beams, similar to \cref{FAxiallyGradedBeamSim1}, is a final test for the patch algorithm.
This example has the initial beam deformation
\begin{align}
u =v =0,\quad
w =0.1x^2\left(3L-x\right) \quad\text{at }t=0\,.
\label{eq:InitDef}
\end{align}
The beam starts vibrating after being released and the simulation period is~$300$ time units.
\cref{fig:DefAlB} displays displacements of the free tip of the beams. 
The $9$~patch and $17$~patch curves are plotted together with the full-domain curve.
The $17$~patch curve closely matches the full-domain curve, whereas the $9$~patch curve has a slightly shorter frequency.
These simulations indicate that the current patch scheme obtains higher accuracy by increasing the number of patches, albeit with increasing computational time.

In \cref{fig:DefAlB} the non-dimensional vibration periods are~$125$ and~$110$ for the three-layer and five-layer beams, respectively. 
These periods correspond to non-dimensional angular frequencies of~$0.05$ and~$0.057$, respectively, which are very close to the frequency of the first $z$-bending mode of the beam listed in \cref{Tab:FreCSGB}, indicating that this mode is the dominant mode in describing the beam dynamics.
With the non-dimensional viscosity coefficient $\eta=10^{-3}$, the vibration amplitude reduces by about~$50\%$ after every period for the three-layer beam, and~$26\%$ for the five-layer beam.
The larger dissipation in the three-layer beam is due to higher stresses between the layers caused by sharper differences in material properties in adjacent layers (recall that volume fractions change by~\(20\%\) from layer to layer, but in the five-layer beam change by~\(10\%\)).


\section{Patch scheme embeds macroscale dynamics}
\label{Spsemd}

The patch scheme does not assume any macroscale variables. 
Consequentially, \emph{a crucial question} is how can we be assured that the patch scheme captures the macroscale slow dynamics and structures?
The \cite{Whitney1936} embedding theorem furnishes an answer.

Roughly, the theorem attests to the fact that every \(m\)D~manifold is parametrisable from almost every subspace of more than~\(2m\)D. 
Here this ensures, in a sense, that the patch scheme provides the higher-D subspace in which the slow manifold of the macroscale viscoelastic beam dynamics is embedded, whatever it may be.  

For 3D beams, the basic macroscale beam models have, at each cross-section,  displacement and velocity of two bending modes as well as torsion and axial compression modes.  
Thus, at each and every cross-section, the 3D viscoelastic beam dynamics has a basic slow manifold that is \(m=8\)D in state space.\footnote{Interpret statements, invoking a manifold or subspace ``at every cross-section", in the sense developed by the theory of \cite{Roberts2013a}.  That is, in systems of large spatial extent there commonly are useful, spatially global, invariant manifolds of high-D that are effectively decomposable into a union of \emph{spatially local} manifolds\slash subspaces of much lower dimension---a dimension determined by the spatial cross-section---and that are weakly coupled to neighbouring locales.}
Alternatively, 3D multi-continua\slash micromorphic\slash micropolar models of beams \cite[e.g.,][]{Forest2011, Somnic2022} include extra assumed cross-sectional modes  to the macroscale modelling.
For example,  micropolar models include three extra variables, the micro-rotations in 3D, leading to six more dynamical dimensions, and hence leading to a quasi-slow manifold of \(m=14\)D in the state space at every cross-section.
Such physically based models are, in effect, approximate slow manifolds because they focus on the relatively slow viscoelastic waves of macroscale wavelength: they neglect all the faster high-frequency cross-waves.

But a significant limitation of the usual approaches to such multi-continua modelling is that the approaches \emph{assume} various cross-sectional structures for the additional modes.
In the presence of functionally graded beams, such as those addressed in \cref{Sec:CrSecGB}, and in order to accurately model the macroscale, the cross-sectional structures of additional modes need to be non-trivially tailored to the particular graded material.
In analytical (asymptotic) homogenisation, multi-mode techniques, such as those invariant manifold models \cite{Watt94b} developed to accurately model shear dispersion, would rigorously and systematically derive algebraic multi-continua homogenisations.
However, such analysis takes considerable effort (current research).
Instead, the patch scheme, aided by the Whitney Embedding Theorem, automatically, implicitly, and painlessly provides the necessary accurate closures.


Consider an example case from the functionally graded beams of \cref{Sec:CrSecGB} where the patch scheme is applied to 3D beams with a cross-section of \((n_y,n_z)=(6,6)\) interior micro-grid points.
Each microscale point has the three displacement variables~\(u,v,w\), so across each patch face there are \(3\cdot n_y\cdot n_z=108\) variables.  
In elasticity, each displacement degree-of-freedom has a corresponding velocity in state space, and so leads to \(216\)~dimensions for the state space of each cross-section.
The patch scheme of \cref{Sec:Patchscheme} communicates the complete information about this \(216\)D space on the left\slash right faces of the patches to realise the macroscale inter-patch coupling.

Because \(216>2\cdot14>2\cdot8\)\,, the Whitney Embedding Theorem asserts that the patch scheme exchanges enough information to almost surely parametrise both such mentioned slow manifolds of the macroscale dynamics (both basic and micropolar).
The patch scheme does \emph{not} need to explicitly compute and exchange  specific assumed macroscale average quantities,  computed via assumed weighting structures.

Consequently, if the microscale system actually has some unknown and/or unsuspected macroscale `average' variables, then Whitney's Embedding Theorem persuades us that the patch scheme generally also captures the macroscale slow dynamics of these unknown and/or unsuspected variables---just provided~\(n_y\cdot n_z\) is not too small.

\section{Axially graded beam}
\label{Sec:AxialGB}

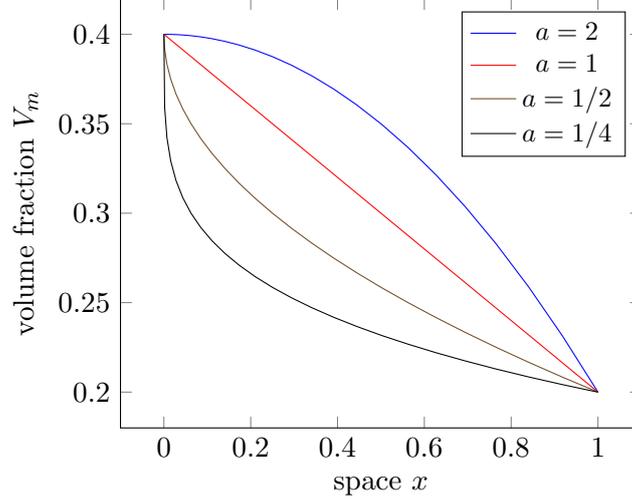
\begin{figure}
\centering
\caption{\label{fig:VFAGB} Volume fraction of~\Al\ along the axis of the axial graded beam in~\cref{eq:VmAGB} for different values of the exponent~$a$.}
\input{Figs/VFAGB.tex}
\end{figure}
An axially graded beam has elastic properties varying 
along the longitudinal \(x\)-axis \citep{Ganta2022}.
In our analysis and example simulations of axially graded beams, we use the same \Al/\Si\C\ properties as those used in \cref{Sec:CrSecGB}.

For beams constructed from~\Al\ and~\Si\C\ we set the volume fraction of the metal component~\Al\ to decrease from~40\% at the fixed end to~20\% at the free end according to
\begin{equation}
V_m(x)=0.2\left(2-x^a\right), \quad 0\le x\le 1\,,
\label{eq:VmAGB}
\end{equation}
in which the exponent~$a$ is a positive real number.
The volume fraction of the ceramic component~\Si\C\ is the complement of the metal component: $V_c(x)=1-V_m(x)$. 
\cref{fig:VFAGB} illustrates the shape of the volume fraction~$V_m(x)$ for different values of exponent~$a$.
Similarly to the case of the cross-sectional graded beam, the Young's modulus of the \Al-\Si\C\ mixture is determined using~\cref{eq:EMixRule}, whereas the mixing rule~\cref{eq:LinMixRule} determines Poison's ratio and the density.
Beam dimensions are~$110\opn{mm}$ in length, and~$10\opn{mm}$ in both width and thickness---in this section the beams have square cross-section. 
As in \cref{Sec:CrSecGB}, a static force of~$147\opn{N}$ applies at the free end of the beams in the static equilibrium computations, and the initial deformation~\cref{eq:InitDef} is used in dynamic simulations of beam vibration such as \cref{FAxiallyGradedBeamSim1}.
Four exponent cases $a=2,1,0.5,0.125$ are investigated.  

\begin{figure}\centering
\caption{\label{fig:DispAGB} Vertical \(z\)-displacement of the axially graded \Al/\Si\C\ beams computed using $9$~patches, $17$~patches, and the full-domain for exponent \(a=2,1,1/2,1/8\) of the volume fraction~\cref{eq:VmAGB} for~\Al.}
\begin{tabular}{@{}l@{}l@{}}
(a) $a=2$ & (b) $a=1$\\
\inPlot{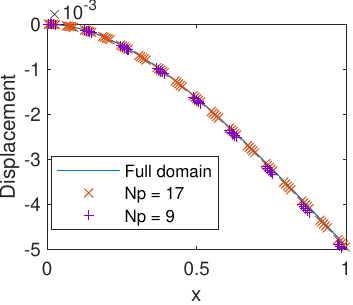}&
\inPlot{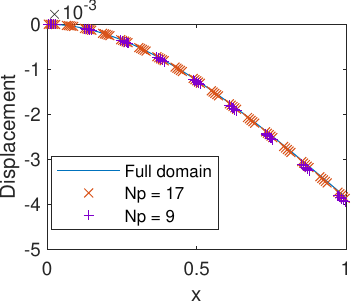}\\
(c) $a=0.5$ & (d) $a=0.125$\\ 
\inPlot{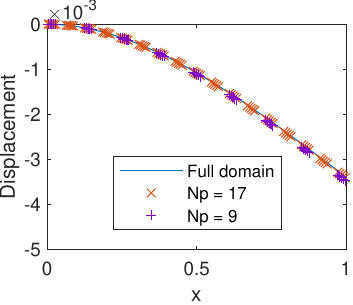}&
\inPlot{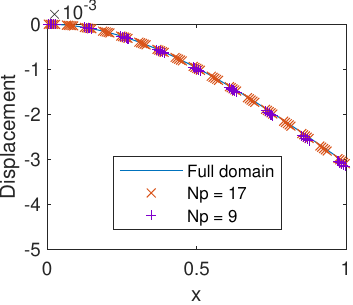}
\end{tabular}
\end{figure}

\begin{figure}\centering
\caption{\label{fig:DispErrorAGB} Relative errors of the vertical displacements  relative to the maximum displacement of the free tip in the corresponding full-domain computations for different values of the exponent \(a=2,1,1/2,1/8\) of volume fraction functions~\cref{eq:VmAGB} for~\Al.}
\begin{tabular}{@{}l@{}l@{}}
(a) $a=2$ & (b) $a=1$\\
\inPlot{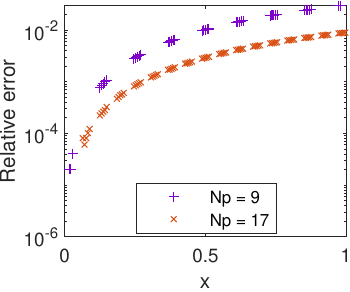}&
\inPlot{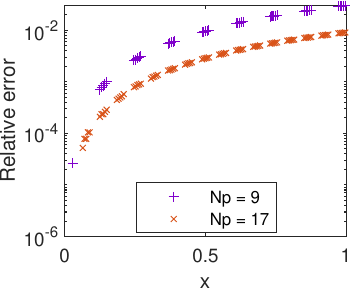}\\
(c) $a=0.5$ & (d) $a=0.125$\\ 
\inPlot{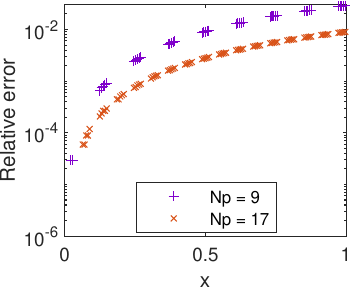}&
\inPlot{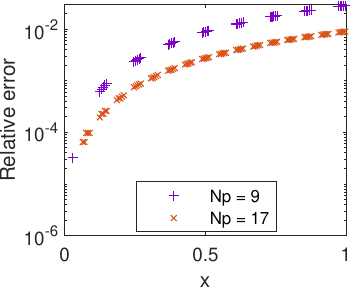}
\end{tabular}
\end{figure}

\cref{fig:DispAGB} plots displacements of the axially graded beam.
The $9$~patch and $17$~patch displacements align with the full-domain displacement in the four cases of exponent~$a$.
The beam deformation decreases with decreasing~$a$, indicating that smaller exponent~$a$ produces a stiffer beam.
Physically this is expected because a smaller~$a$ more rapidly decreases the volume fraction of~\Al\ (\cref{fig:VFAGB}), and as the Young's modulus of~\Si\C\ is about~$4.5$ times the Young's modulus of~\Al, a beam with a larger volume fraction of~\Si\C\ will be more resistant to deformation when under load.  
In the $9$~patch computations, the maximum displacements at the tip of the free end are $-5.06\xten{-3}$, $-3.95\xten{-3}$, $-3.48\xten{-3}$, and $-3.16\xten{-3}$ for $a=2$, $1$, $0.5$, and~$0.125$, respectively. 
In the $17$~patch computations, the maximum displacements are $-4.95\xten{-3}$, $-3.87\xten{-3}$, $-3.41\xten{-3}$, and~$-3.10\xten{-3}$, respectively.   
\cref{fig:DispErrorAGB} plots errors of beam displacements, calculated relative to the maximum displacement of the corresponding full-domain computation.
In all the cases of exponent~$a$, the errors are less than~$3\%$ with $9$~patches, and less than~$1\%$ with $17$~patches.

\begin{figure}\centering
\caption{\label{fig:EValueAGB17} Leading~4\% of the eigenvalues of the  axial graded beams for different power indices of the volume fraction functions of~\Al\ computed using $17$~patches and exponent \(a=2,1,1/2,1/8\) of volume fraction~\cref{eq:VmAGB}.  
The~$\circ$, {\tiny$\square$}, $\triangleright$, and~{\small$\vartriangle$} symbols refer to the $z$-bending modes, $y$-bending modes, torsion modes (in the $yz$-plane), and axial compression modes, respectively: red, green, and yellow colours denote the first, second and third modes.}
\begin{tabular}{@{}l@{}l@{}}
(a) $a=2$ & (b) $a=1$\\
\inPlot{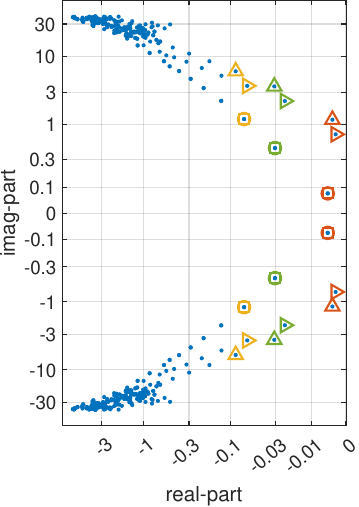}&
\inPlot{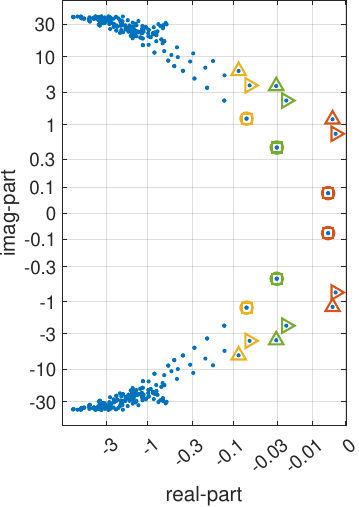}\\
(c) $a=0.5$ & (d) $a=0.125$\\ 
\inPlot{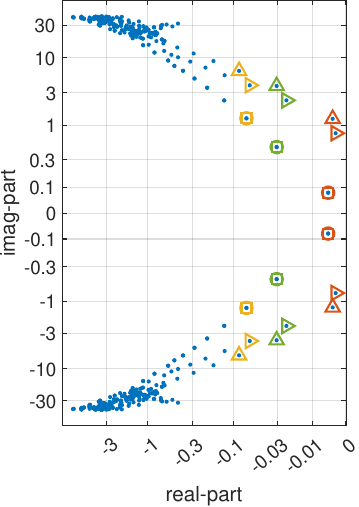}&
\inPlot{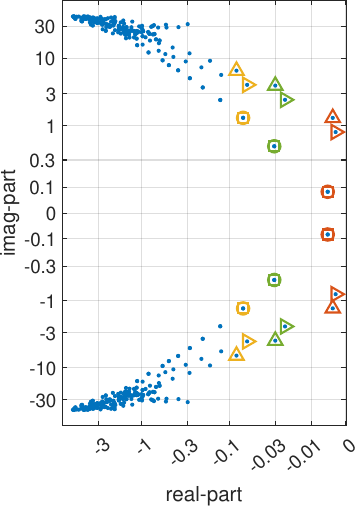}\\ 
\end{tabular}
\end{figure}

\cref{fig:EValueAGB17} plots eigenvalues of the numerically computed Jacobian matrix of the coded equations for the \(17\)~patch scheme applied to the axially graded beam. 
Recall that for every mode, the real part of the eigenvalues reflects the viscoelastic dissipation, and the imaginary-part is the frequency.    
We again observe modes of~$y$- and $z$-bending, torsion, and axial compression.
\cref{Tab:EigvalAGB} lists the frequencies of the first three modes of these four vibration types.
Computed eigenvalues of the $y$-bending and $z$-bending are equal, as is expected because  the cross-section of the beam is square. 
The frequencies of all vibration types and modes consistently increase as the exponent~$a$ decreases, which is expected because ceramic beams typically have higher frequencies of oscillation than metal beams.  
Frequencies from the $9$~patch and $17$~patch computations agree well, with small differences increasing slightly with higher frequencies.
\begin{table}\centering
\caption{\label{Tab:EigvalAGB}Eigenvalues give the following frequencies for the first three modes of the $z$-bending, $y$-bending, torsion, and axial compression.  
These frequencies are for exponent \(a=2,1,1/2,1/8\) of the volume fraction~\cref{eq:VmAGB}.}
\begin{tabular}{lcccccc}
\hline
& & \multicolumn{4}{c}{Angular frequency $\omega^*$}\\
         &  & $a=2$ & $a=1$ & $a=0.5$ & $a=0.125$ \\
\hline
         &  & \multicolumn{4}{c}{9-patch computation} \\
\hline  
\multirow{3}{*}{$y$- and $z$-bending} & 1 & $0.073$ & $0.073$ & $0.077$ & $0.081$\\
            & 2 & $0.493$ & $0.486$ & $0.486$ & $0.498$\\
            & 3 & $0.871$ & $0.945$ & $0.988$ & $1.025$\\          
\hline
            & 1 & $0.716$ & $0.733$ & $0.758$ & $0.800$\\
torsion     & 2 & $2.233$ & $2.285$ & $2.342$ & $2.430$\\
            & 3 & $3.935$ & $3.965$ & $4.015$ & $ 4.117$\\          
\hline
            & 1 & $1.176$ & $1.203$ & $1.240$ & $1.302$\\
compression & 2 & $3.660$ & $3.737$ & $3.822$ & $3.953$\\
            & 3 & $6.420$ & $6.469$ & $ 6.543$ & $ 6.692$\\          
\hline
         &  & \multicolumn{4}{c}{17-patch computation} \\  
\hline
\multirow{3}{*}{$y$- and $z$-bending} & 1 & $0.072$ & $0.072$ & $0.075$ & $0.080$\\
 & 2 & $0.446$ & $0.446$ & $0.469$ & $0.488$ \\
            & 3 & $1.203$ & $1.203$ & $1.260$ & $1.305$\\         
      
\hline
            & 1 & $0.716$ & $0.716$ & $0.758$ & $0.800$\\
torsion     & 2 & $2.216$ & $2.216$ & $2.332$ & $2.420$\\
            & 3 & $3.706$ & $3.706$ & $3.899$ & $4.041$\\          
\hline
            & 1 & $1.176$ & $1.176$ & $1.239$ & $1.302$\\
compression & 2 & $3.631$ & $3.631$ & $3.804$ & $3.935$\\
            & 3 & $6.062$ & $6.062$ & $6.351$ & $6.562$\\          
\hline
\end{tabular}
\end{table}

\begin{figure}\centering
\caption{\label{fig:TVAGB} Displacement at the free end of the axially graded~\Al/\Si\C\ beam simulated using $9$~patches, $17$~patches, and the full-domain for exponents \(a=2,1,1/2,1/8\) of the volume fraction~\cref{eq:VmAGB} for~\Al:
the blue solid curves are the full-domain simulations;  symbol $\color{matlab2}\times$ for $17$~patch simulations; symbol $\color{matlab1}+$ for $9$~patch simulations.}
\begin{tabular}{@{}l@{}l@{}}
(a) $a=2$ & (b) $a=1$\\
\inPlot{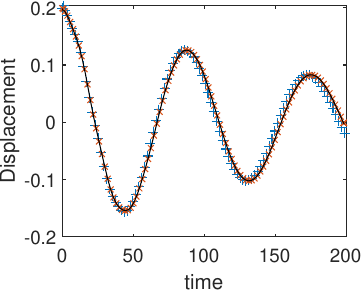}&
\inPlot{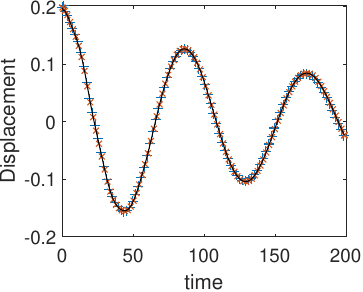}\\
(c) $a=0.5$ & (d) $a=0.125$\\ 
\inPlot{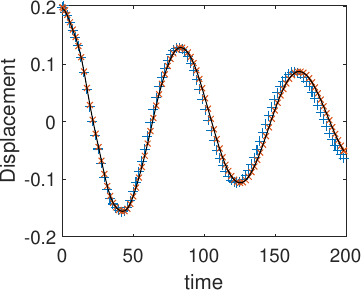}&
\inPlot{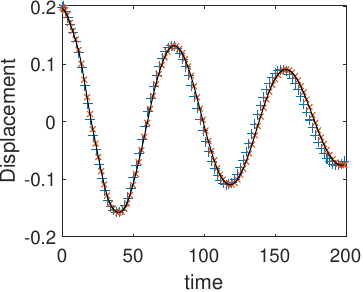}  
\end{tabular}
\end{figure}

\cref{fig:TVAGB} plots displacements of the free end versus time obtained from simulations of vibrating beams with initial deformation~\cref{eq:InitDef}.
In all cases of exponent~$a$, the $9$~patches, $17$~patches, and full-domain curves are very similar.
The vibration period (non-dimensional) of the free end is~$94$, $85$, $80$, and~$77$ for $a=2$, $1$, $0.5$, and~$0.125$, respectively.
Equivalently, the vibrating frequency is~$0.067$, $0.074$, $0.078$, and~$0.082$, respectively, which is the frequency of the first $z$-bending modes, which is also the lowest frequency mode in the beam. 
As expected, as exponent~$a$ is reduced, the beam vibrates with a higher frequency.
Viscoelastic dissipation is evident in the decay over time of the vibration amplitude.   
The vibration amplitude reduces about~$41\%$, $37\%$, $34\%$ and~$29\%$, respectively, over each period of vibration. 

When manufacturing an axially graded beam, the volume fractions will not vary as smoothly as described by \cref{eq:VmAGB}.
Instead, the manufactured volume fractions will randomly vary about a designed value.
To account for such random variations, we introduced a random term into the Young's modulus:
\begin{equation}
E_{\textsc{ragb}}(x,y,z):=E_{\textsc{agb}}(x)\left(1+\alpha U[-1,1]\right),\label{eq: randFElas}
\end{equation}  
in which $U[-1,1]$ is a random number chosen uniformly in the interval~$[-1,1]$ for each \((y,z)\)~in the cross-section, $\alpha$~is the amplitude of the random variation, and $E_{\textsc{agb}}(x)$ is the functional Young's modulus used in the previous axially graded beam simulations.
That is, at each~\(x\) the Young's modulus is uniformly distributed between $E_{\textsc{agb}}(x)\left(1-\alpha\right)$ and $E_{\textsc{agb}}(x)\left(1+\alpha\right)$.
The simulation of \cref{FAxiallyGradedBeamSim1} is for such a randomly perturbed axially graded beam.
 
\cref{Tab:EigvalAGBwRF} lists non-dimensional angular frequencies of the first three modes of the four vibration types.
Firstly, the frequencies of the $y$-bending and $z$-bending vibrations are no longer identical, since the random variations break the \(yz\)-symmetry.
Secondly, differences in the first and second modes of all vibration types are relatively small, whereas there are significant differences in the third modes. 
These differences are due to high frequency vibrations being more sensitive to local heterogeneities in the material.
Finally, frequencies computed using $9$~patches are very consistent with those computed using $17$~patches.  

\begin{figure}
\centering
\caption{\label{fig:TVAGBwRF} Displacement at the free end of the axially graded~\Al/\Si\C\ beam with randomly fluctuated elasticity \cref{eq: randFElas} simulated using $9$~patches, $17$~patches, and  the full-domain:
the blue solid curves are the full-domain simulations; symbol $\color{matlab2}\times$ for $17$~patch simulations; symbol $\color{matlab1}+$ for $9$~patch simulations.}
\begin{tabular}{@{}l@{}l@{}}
(a) $a=2$ & (b) $a=1$\\
\inPlot{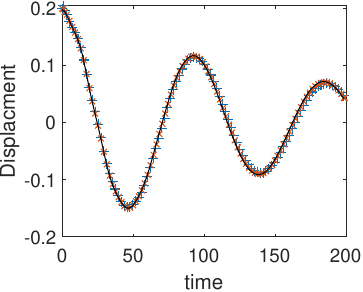}&
\inPlot{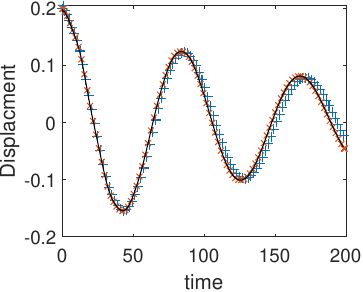}\\
(c) $a=0.5$ & (d) $a=0.125$\\ 
\inPlot{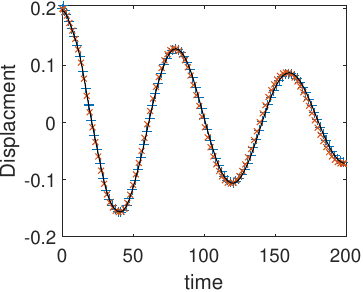}&
\inPlot{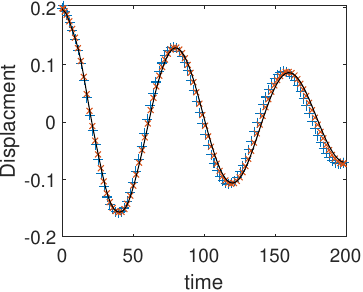}   
\end{tabular}
\end{figure}

\cref{fig:TVAGBwRF}  plots displacements of the free end versus time obtained from simulating the dynamics of axially graded beams with random variations from initial deformation~\cref{eq:InitDef}.
The vibration period of the four cases of \(a=2,1,1/2,1/4\) are $92$, $84$, $79$ and $78$, respectively.  
These values are only slightly different from the corresponding values in the non-random cases (\cref{fig:TVAGB}).     
These slight differences are not surprising as the random fluctuations in the Young's modulus only significantly effect the high frequency modes (\cref{Tab:EigvalAGBwRF}), whereas it is the lowest frequency modes that dominate the dynamics from this initial condition. 
In conclusion, the effects of such random heterogeneity in the beam is straightforwardly predicted \text{with patch dynamics.}
      
\begin{table}
\caption{\label{Tab:EigvalAGBwRF}Eigenvalues give the following frequencies for the first three modes of the $z$-bending , $y$-bending, torsion, and compression modes of vibrations of axially graded beams with the randomly varying elasticity~\cref{eq: randFElas}.  
These frequencies are for exponent \(a=2,1,1/2,1/8\) of the volume fraction~\cref{eq:VmAGB}.}
\centering
\begin{tabular}{lcccccc}
\hline
& & \multicolumn{4}{c}{Angular frequency $\omega^*$}\\
         &  & $a=2$ & $a=1$ & $a=0.5$ & $a=0.125$ \\
\hline
         &  & \multicolumn{4}{c}{9 patch computation} \\
\hline  
            & 1 & $0.068$ & $0.073$ & $0.078$ & $0.079$\\
$z$-bending & 2 & $0.494$ & $0.492$ & $0.502$ & $0.488$\\
            & 3 & $1.006$ & $0.997$ & $1.018$ & $0.987$\\          
\hline
            & 1 & $0.069$ & $0.075$ & $0.080$ & $0.081$\\
$y$-bending & 2 & $0.502$ & $0.505$ & $0.511$ & $0.502$\\
            & 3 & $1.024$ & $1.031$ & $1.049$ & $1.029$\\          
\hline

            & 1 & $0.698$ & $0.747$ & $0.784$ & $0.801$\\
torsion     & 2 & $2.401$ & $2.403$ & $2.441$ & $2.431$\\
            & 3 & $4.122$ & $4.085$ & $4.129$ & $4.098$\\          
\hline
            & 1 & $1.191$ & $1.254$ & $1.302$ & $1.308$\\
compression & 2 & $4.030$ & $4.002$ & $4.035$ & $3.967$\\
            & 3 & $6.905$ & $6.799$ & $6.824$ & $6.687$\\          
\hline
         &  & \multicolumn{4}{c}{17 patch computation} \\  
\hline
            & 1 & $0.068$ & $0.075$ & $0.079$ & $0.078$\\
$z$-bending & 2 & $0.471$ & $0.484$ & $0.495$ & $0477$\\
            & 3 & $1.306$ & $1.314$ & $1.331$ & $1.278$\\          
\hline
            & 1 & $0.069$ & $0.076$ & $0.080$ & $0.080$\\
$y$-bending & 2 & $0.477$ & $0.491$ & $0.502$ & $0.489$\\
            & 3 & $1.322$ & $1.334$ & $1.348$ & $1.308$\\          
\hline

            & 1 & $0.702$ & $0.753$ & $0.786$ & $0.820$\\
torsion     & 2 & $2.402$ & $2.412$ & $2.436$ & $2.481$\\
            & 3 & $4.040$ & $4.041$ & $4.073$ & $4.139$\\          
\hline
            & 1 & $1.203$ & $1.273$ & $1.329$ & $1.287$\\
compression & 2 & $4.052$ & $4.043$ & $4.102$ & $3.886$\\
            & 3 & $6.801$ & $6.762$ & $6.848$ & $6.475$\\          
\hline
\end{tabular}
\end{table}

\section{Conclusion}
\label{Sec:Conclusion}

Herein we further developed a novel multiscale patch dynamics algorithm to investigate both static and dynamic characteristics of cross-sectionally and axially graded cantilever beams made of~\Al/\Si\C. 
The patch dynamics developed here extends our previous application of a patch scheme to two-dimensional heterogeneous beams \citep{Tran2024}.
 
The developed multiscale patch scheme accurately predicts static deflections, vibration modes, and time evolution of functional graded beams made from~\Al/\Si\C\ mixtures.
The computational efficiency of the algorithm is observed when simulating only~$52\%$ and~$27\%$ of the beam's full domain, but still correctly predicting the beam's deformations and vibrations.
When simulating over~$27\%$ of the beam's domain, the patch algorithm is roughly seventeen times faster than a full-domain computation.
Numerical stability is another advantage of the patch scheme, ensured via the carefully designed and proven patch coupling.   
A toolbox to implement the patch scheme is freely available \citep{Roberts2019b}.  

Lastly, the proposed patch scheme makes very few assumptions regarding the form of the macroscale solution, in contrast to other beam homogenisation methods.

\paragraph{Acknowledgements}
This research was supported by Australian Research Council grants DP220103156 and DP200103097.


\end{document}

%% file: preamble.tex
\IfFileExists{ajr.sty}{
\usepackage[a5paper,margin=10mm,tmargin=15mm]{geometry} 
}{
}
\usepackage[collision]{chemsym}

\usepackage{amsmath,amssymb,microtype}
\allowdisplaybreaks
\usepackage[defaultlines=3,all]{nowidow}
\usepackage[leftcaption]{sidecap}
\usepackage{multirow}
\usepackage{xcolor}
\usepackage{natbib}
\AtEndDocument{\bibliography{bibexport}}

\usepackage{tikz,pgfplots}
\pgfplotsset{compat=newest}
\def\tikzsetnextfilename#1{} 
\usepgfplotslibrary{external}

\usepackage{doi}
\hypersetup{
    colorlinks   = true, 
    urlcolor     = blue, 
    linkcolor    = blue, 
    citecolor   = magenta,
    } 
\usepackage[nameinlink,capitalise,noabbrev]{cleveref}
\crefname{equation}{}{}
\let\ref\cref
\let\eqref\cref
\let\autoref\cref

\usepackage[colorinlistoftodos,textsize=footnotesize]{todonotes}
  \setuptodonotes{inline,color=orange!20}
  \ifx\undefined
     
  \else
     
  \fi

\def\RaisedName#1{}
\def\ajrSplit#1#2\ajrEndSplit{\def\ajrcar{#1}\def\ajrcdr{#2}}
\def\Vec{}
\renewcommand{\Vec}[1]{%
  \typeout{**** Command #1v denotes vector #1}%
  \expandafter\ajrSplit#1\ajrEndSplit%
  \ifx\ajrcdr\empty 
    \expandafter\DeclareRobustCommand\csname#1v\endcsname%
    {{\RaisedName{\ensuremath{\backslash}#1v}
      \ensuremath{\ifx#1i \vec\imath
        \else\ifx#1j \vec\jmath
        \else\vec #1\fi\fi}
      \global\csname UsedVec#1true\endcsname
    }}%
  \else 
    \expandafter\DeclareRobustCommand\csname#1v\endcsname%
    {{\ensuremath{\vec{\csname#1\endcsname}}}
      \global\csname UsedVec#1true\endcsname}%
  \fi%
  \expandafter\newif\csname ifUsedVec#1\endcsname
  \csname UsedVec#1false\endcsname
  \AtEndDocument{\csname ifUsedVec#1\endcsname\else
  \typeout{**** Vec symbol #1v not used.}
  \fi}
  }
\newcommand{\Cal}[1]{
  \typeout{**** Command c#1 denotes calligraphic #1}%
  \expandafter\DeclareRobustCommand\csname c#1\endcsname%
  {{\RaisedName{\detokenize{\c}\hspace{-0.5em}#1}
  \ensuremath{\mathcal #1}
  \global\csname UsedCal#1true\endcsname  
  }}
  \expandafter\newif\csname ifUsedCal#1\endcsname
  \csname UsedCal#1false\endcsname
  \AtEndDocument{\csname ifUsedCal#1\endcsname\else
  \typeout{**** Cal symbol c#1 not used.}
  \fi}
  }

\renewcommand{\vec}[1]{\text{\boldmath$#1$}}
\newcommand{\dx}{\delta_x}
\newcommand{\dy}{\delta_y}
\newcommand{\dz}{\delta_z}
\newcommand{\E}{\cdot 10^}

\def\matlab{\textsc{Matlab}}

\def\pde{\textsc{pde}}

\def\rve{\textsc{rve}}
\Vec e\Vec f\Vec k\Vec n\Vec q\Vec x\Vec u
\Vec{varepsilon}\Vec{sigma}\Vec I\Vec J\Vec F
\Cal U

\newcommand{\Dn}[3]{\mathchoice
  {\frac{\partial^{#2} #3}{\partial #1^{#2}}}%
  {{\partial^{#2} #3}/{\partial #1^{#2}}}%
  {{\partial^{#2} #3}/{\partial #1^{#2}}}%
  {{\partial^{#2} #3}/{\partial #1^{#2}}}%
  }
\newcommand{\DD}[2]{\Dn{#1}2{#2}}

\newcommand{\Ord}[1]{\ensuremath{\mathcal O%
  \mathchoice{\big(#1\big)}{\big(#1\big)}{(#1)}{(#1)}%
  }}

\newcommand{\opn}[1]{\operatorname{#1}}

\def\E#1{\textsc{e}\ifx#1-{-}\else\ifx#1+{+}\else#1\fi\fi}
\def\xten{\cdot10^}

\def\uc{blue}
\def\vc{red}
\def\wc{green!60!black}
\def\sc{black}

\def\vLen{0.2} 
\def\o{0.03}   

\def\phalf{^+} 
\def\mhalf{^-} 
\def\kbk{\textsc{kbk08}}

\definecolor{matlab1}{rgb}{0,0.4470,0.7410}
\definecolor{matlab2}{rgb}{0.8500,0.3250,0.0980}
\definecolor{matlab3}{rgb}{0.9290,0.6940,0.1250}
\definecolor{matlab4}{rgb}{0.4940,0.1840,0.5560}
\definecolor{matlab5}{rgb}{0.4660,0.6740,0.1880}
\definecolor{matlab6}{rgb}{0.3010,0.7450,0.9330}
\definecolor{matlab7}{rgb}{0.6350,0.0780,0.1840}

%% file: Figs/3DCell.tex

\begin{tikzpicture}
\tikzset{every mark/.append style={scale=2}}
\begin{axis}[axis equal image,view={20}{30}
    ,xlabel={$x$},ylabel={$y$},zlabel={$z$}
    ,xtick=\empty,ytick=\empty,ztick=\empty
    ,axis line style={draw=none}
    ,xmin=-0.1,ymin=-0.1,zmin=-0.1
    ,xmax=1.3,ymax=2.1,zmax=1.5]
     \draw (0,0,0) -- (0,2,0) -- (1.2,2,0) -- (1.2,0,0) -- cycle;
     \draw (0,0,1.4) -- (0,2,1.4) -- (1.2,2,1.4) -- (1.2,0,1.4) -- cycle;
     \draw (0,0,0) -- (0,0,1.4);
     \draw (0,2,0) -- (0,2,1.4);
     \draw (1.2,0,0) -- (1.2,0,1.4);
     \draw (1.2,2,0) -- (1.2,2,1.4);    
\foreach \x in {0,1}{ 
\addplot3[mark=oplus,\sc] coordinates {(1.2*\x,1,0.7)};
\addplot3[mark=halfcircle,\sc] coordinates {(0.5,2*\x,0.7)};
\addplot3[mark=o,\sc] coordinates {(0.6,1.0,1.4*\x)};
\foreach \y in {0,1}{
\addplot3 [quiver={u=0,v=0,w=\vLen},\wc, -stealth,ultra thick]
          coordinates {(1.2*\x,2*\y,0.7-\o)};
\addplot3 [quiver={u=0,v=\vLen,w=0},\vc, -stealth,ultra thick]
          coordinates {(1.2*\x,1,1.4*\y)};
\addplot3 [quiver={u=\vLen,v=0,w=0},\uc, -stealth,ultra thick]
          coordinates {(0.6,2*\x,1.4*\y)}; 
\foreach \z in {0,1}{
\addplot3[mark=otimes,\sc] coordinates {(1.2*\x,2*\y,1.4*\z)};
}}}
\end{axis}
\end{tikzpicture}

%% file: Figs/xyplane.tex
\begin{tikzpicture}
\begin{axis}[axis equal image
    ,xlabel={$x$},ylabel={$y$}
    ,xtick=\empty,ytick=\empty
    ,axis line style={draw=none}
    ,xmin=-1.2,ymin=-1.2
    ,xmax=7.3,ymax=7.3]
    \draw (-1,0) -- (-1,6) -- (5.5,6) -- (5.5,0) -- cycle;
\foreach \x in {0,...,5}{
\addplot [mark=square,teal] coordinates {(\x,-0.6)};
\addplot [mark=square,teal] coordinates {(\x,6.6)};
\addplot [mark=square,teal] coordinates {(\x-0.5,-0.6)};
\addplot [mark=square,teal] coordinates {(\x-0.5,6.6)};
\addplot [mark=square,orange] coordinates {(\x,0)};
\addplot [mark=square,orange] coordinates {(\x,6)};
\addplot [mark=square,black] coordinates {(-0.5,1.2*\x)};
\addplot [mark=square,black] coordinates {(5.5,1.2*\x)};
\foreach \y in {0,...,6}{ 
\addplot [quiver={u=0,v=\vLen},\vc, -stealth,ultra thick]
          coordinates {(\x,1.2*\y-0.6)};
\addplot [quiver={u=\vLen,v=0},\uc, -stealth,ultra thick]
          coordinates {(\y-0.5,1.2*\x)};         
\addplot[mark=o,\sc] coordinates {(\x-0.5,1.2*\y-0.6)};          
}
\foreach \y in {0,...,5}{ 
\addplot[mark=otimes,\sc] coordinates {(\x,1.2*\y)};
}}
\foreach \y in {0,...,4}{
\addplot [quiver={u=0,v=\vLen},\vc, -stealth,ultra thick]
          coordinates {(-1,1.2*\y+0.6)};
\addplot [mark=square,black] coordinates {(-1,1.2*\y+0.6)};
\addplot [mark=square,black] coordinates {(5,1.2*\y+0.6)};
}
\end{axis}
\end{tikzpicture}

%% file: Figs/xzplane.tex
\begin{tikzpicture}
\begin{axis}[axis equal image
    ,xlabel={$x$},ylabel={$z$}
    ,xtick=\empty,ytick=\empty
    ,axis line style={draw=none}
    ,xmin=-1.2,ymin=-1.2
    ,xmax=7.3,ymax=7.3]
    \draw (-1,0) -- (-1,5) -- (5.5,5) -- (5.5,0) -- cycle;
\foreach \x in {0,...,5}{
\addplot [mark=square,teal] coordinates {(\x,-0.5)};
\addplot [mark=square,teal] coordinates {(\x,5.5)};
\addplot [mark=square,teal] coordinates {(\x-0.5,-0.5)};
\addplot [mark=square,teal] coordinates {(\x-0.5,5.5)};
\addplot [mark=square,orange] coordinates {(\x,0)};
\addplot [mark=square,orange] coordinates {(\x,5)};
\addplot [mark=square,black] coordinates {(-0.5,\x)};
\addplot [mark=square,black] coordinates {(5.5,\x)};
\foreach \y in {0,...,6}{ 
\addplot [quiver={u=0,v=\vLen},\wc, -stealth,ultra thick]
          coordinates {(\x,\y-0.5)};
\addplot [quiver={u=\vLen,v=0},\uc, -stealth,ultra thick]
          coordinates {(\y-0.5,\x)};         
\addplot[mark=halfcircle,\sc] coordinates {(\x-0.5,\y-0.5)};          
}
\foreach \y in {0,...,5}{ 
\addplot[mark=otimes,\sc] coordinates {(\x,\y)};
}}
\foreach \y in {0,...,4}{
\addplot [quiver={u=0,v=\vLen},\wc, -stealth,ultra thick]
          coordinates {(-1,\y+0.5)};
\addplot [mark=square,black] coordinates {(-1,\y+0.5)};
\addplot [mark=square,black] coordinates {(5,\y+0.5)};
}
\end{axis}
\end{tikzpicture}

%% file: Figs/yzplane.tex
\begin{tikzpicture}
\begin{axis}[axis equal image
    ,xlabel={$y$},ylabel={$z$}
    ,xtick=\empty,ytick=\empty
    ,axis line style={draw=none}
    ,xmin=-1.2,ymin=-1.2
    ,xmax=7.3,ymax=7.3]
    \draw (-0.6,0) -- (-0.6,5) -- (5.4,5) -- (5.4,0) -- cycle;
\foreach \x in {0,...,5}{
\foreach \y in {0,...,6}{ 
\addplot [quiver={u=0,v=\vLen},\wc, -stealth,ultra thick]
          coordinates {(1.2*\x-0.6,\y-0.5)};
\addplot [quiver={u=\vLen,v=0},\vc, -stealth,ultra thick]
          coordinates {(1.2*\y-1.2,\x)};              
}
\foreach \y in {0,...,5}{
\addplot[mark=otimes,\sc] coordinates {(1.2*\x-0.6,\y)};         
}}
\foreach \x in {0,...,4}{
\addplot[mark=oplus,\sc] coordinates {(1.2*\x,-0.5)};  
\addplot[mark=oplus,\sc] coordinates {(1.2*\x,5.5)};
\addplot[mark=oplus,\sc] coordinates {(-1.2,\x+0.5)};
\addplot[mark=oplus,\sc] coordinates {(6,\x+0.5)};
\addplot[mark=square,teal] coordinates {(1.2*\x,-0.5)};
\addplot[mark=square,teal] coordinates {(1.2*\x,5.5)};
\addplot[mark=square,teal] coordinates {(-1.2,\x+0.5)};
\addplot[mark=square,teal] coordinates {(6,\x+0.5)};
\foreach \y in {0,...,4}{
\addplot[mark=oplus,\sc] coordinates {(1.2*\x,\y+0.5)};  
}}
\foreach \x in {0,...,5}{
\addplot[mark=square,teal] coordinates {(1.2*\x-0.6,-0.5)};
\addplot[mark=square,teal] coordinates {(1.2*\x-0.6,5.5)};
\addplot[mark=square,teal] coordinates {(-1.2,\x)};
\addplot[mark=square,teal] coordinates {(6,\x)};
\addplot[mark=square,orange] coordinates {(-0.6,\x)};
\addplot[mark=square,orange] coordinates {(5.4,\x)};
}
\foreach \x in {0,...,3}{
\addplot[mark=square,orange] coordinates {(1.2*\x+0.6,0)};
\addplot[mark=square,orange] coordinates {(1.2*\x+0.6,5)};
}
\end{axis}
\end{tikzpicture}

%% file: Figs/figpatchscheme.tex
\tikzsetnextfilename{Figs/figpatchscheme}
\def\N{7} 
\begin{tikzpicture}[x=18mm,y=18mm]
\foreach \j in {1,...,\N} {
       \filldraw[fill=green!15,draw=green,thick] (\j-0.25,-0.1)--(\j+0.25,-0.1)--(\j+0.25,0.1)--(\j-0.25,0.1)--(\j-0.25,-0.1);   
       \draw (\j-0.32,-0.1) rectangle (\j-0.28,0.1);
       \draw (\j+0.32,-0.1) rectangle (\j+0.28,0.1);
       }
\foreach \j in {1,2,3} {
       \node[] at (4-\j,0) {$I-\j$};  
       \node[] at (4+\j,0) {$I+\j$};  
      }   
\node[] at (4,0) {$I$};           
\draw [brown, thin] (0.5,0.1)--(\N+0.5,0.1);
\draw [brown, thin] (0.5,-0.1)--(\N+0.5,-0.1);
\draw[black,<->] (1,-0.2)-- node[midway, below] {$H$} (2,-0.2);    
\draw[black,<->] (6-0.25,-0.2)-- node[midway, below] {$h$} (6+0.25,-0.2);    
{\color{red}
\draw [->] (2+0.2,0.1) to [out=30,in=90] (4-0.3,0.1);  
\draw [->] (3+0.2,0.1) to [out=30,in=90] (4-0.3,0.1);
\draw [->] (4+0.2,0.1) to [out=150,in=90] (4-0.3,0.1);
\draw [->] (5+0.2,0.1) to [out=150,in=90] (4-0.3,0.1);
\draw [->] (6+0.2,0.1) to [out=150,in=90] (4-0.3,0.1);
}
{\color{blue}
\draw [->] (2-0.2,-0.1) to [out=-30,in=-90] (4+0.3,-0.1);
\draw [->] (3-0.2,-0.1) to [out=-30,in=-90] (4+0.3,-0.1);
\draw [->] (4-0.2,-0.1) to [out=-30,in=-90] (4+0.3,-0.1);
\draw [->] (5-0.2,-0.1) to [out=-150,in=-90] (4+0.3,-0.1);
\draw [->] (6-0.2,-0.1) to [out=-150,in=-90] (4+0.3,-0.1);
}
\end{tikzpicture}  

%% file: Figs/VFAGB.tex
\begin{tikzpicture}
\begin{axis}[no marks,
    xlabel={space $x$},ylabel={volume fraction $V_m$},
    legend pos=north east,legend style={font=\small},domain=0:1]
\foreach \a in {2,1,1/2,1/4} {
    \addplot+[] ({x^2},{0.2*(2-x^(2*\a))});
    \addlegendentryexpanded{$a=\a$}
}
\end{axis}
\end{tikzpicture}